\documentclass[10pt,final,twocolumn]{IEEEtran}
\pdfoutput=1
\usepackage{amssymb,amsmath}
\newcommand{\bea}		{\begin{eqnarray}}
\newcommand{\eea}		{\end{eqnarray}}
\newcommand{\be}		{\begin{equation}}
\newcommand{\ee}		{\end{equation}}
\newcommand{\ba}		{\begin{array}}
\newcommand{\ea}		{\end{array}}

\newcommand{\Ltwo}		{{\mathcal{L}_{2}}}
\newcommand{\ltwo}		{{$\Ltwo$}}
\newcommand{\Ltwoloc}	{{\mathcal{L}_2^e}}
\newcommand{\ltwoloc}	{{$\Ltwoloc$}}
\def\g{\gamma}
\def\k{\beta}

\usepackage{graphicx}
\usepackage{amssymb,amsmath}
\usepackage{epstopdf}

\def\K{\mathcal{K}}
\def\Kinf{\K_{\infty}}
\def\KL{\mathcal{KL}}
\def\L{\mathcal{L}}

\def\R{\mathbb{R}}

\def\Rnn{\R_{\geq 0}}

\def\W{\mathcal{W}}

\def\Z{\mathbb{Z}}

\newcommand{\halmos}{\hfill \blacksquare}
\newcommand{\nn}{{\nonumber}}

\newcommand{\twonorm}[1]{|\!|#1|\!|_{\L_2[0,t]}^2}

\newtheorem{lemma}{Lemma}
\newtheorem{thm}{Theorem}
\newtheorem{prop}{Proposition}

\newtheorem{defn}{Definition}

\newtheorem{remark}{Remark}

\newtheorem{example}{Example}

\title{Input-to-State Stability, integral Input-to-State Stability, and \ltwo-Gain Properties: Qualitative Equivalences
and Interconnected Systems\thanks{
Preliminary versions of this paper appeared at the {\it 51st IEEE Conference on Decision and Control}
and the 2013 {\it Australian Control Conference}.
}
}

\author{Christopher M. Kellett\thanks{C.M.~Kellett is with the School of Electrical Engineering and Computer
Science, University of Newcastle, Callaghan, New South Wales, Australia {\tt\small Chris.Kellett@newcastle.edu.au}.
Kellett is supported by ARC Future Fellowship FT1101000746 and by the Alexander von Humboldt
Foundation.}  \ and
Peter M. Dower\thanks{
 P.M.~Dower 
is with the Department of Electrical \& Electronic Engineering, University of Melbourne, Melbourne, Victoria, Australia 
{\tt\small pdower@unimelb.edu.au}.  Dower is supported by AFOSR grant FA2386-12-1-4084
and ARC Discovery Project 120101549. }
}

\begin{document}
\maketitle

\begin{abstract}
	Input-to-state stability (ISS) and $\L_2$-gain are well-known robust stability properties that
	continue to find wide application in the analysis and control of nonlinear dynamical systems and their
	interconnections.
	We investigate the relationship between ISS-type
	 and $\L_2$-gain properties, demonstrating several qualitative equivalences
	between these two approaches.  We subsequently present several new sufficient conditions for the 
	stability of interconnected systems derived by exploiting these qualitative equivalences.
\end{abstract}

\section{Introduction}
Historically, there have been two dominant approaches to the study of interconnected dynamical systems via a modular, stability-of-subsystems approach.  The first was pioneered by Zames in the 1960's and employs $\L_2$-gain from the input to the output of a (sub-)system \cite{Zame66-TAC} (see also \cite{DeVi75}).  The second approach developed from the introduction of the Input-to-State Stability (ISS) concept by Sontag in 1989 \cite{Sont89-TAC}
which extended classical state-space stability notions for systems described by ordinary differential
equations to include inputs.

The $\L_2$-gain input-output approach, derived largely from frequency domain considerations, led
to the highly successful linear $\mathcal{H}^\infty$ optimal control techniques for linear 
time-invariant systems first suggested in \cite{Zame81-TAC} where it is possible to design feedback 
controllers in a systematic way to achieve a desired closed-loop $\L_2$-gain from disturbance input to 
some suitably weighted penalty variable (regarded as an output) (see \cite{ZDG96} and the references
therein).  The state space formulation and
solution of the $\mathcal{H}^\infty$ optimal control problem for linear systems \cite{DGKF89-TAC}
paved the way for extending these techniques to the study of nonlinear systems where the design
goal remained the design of a closed-loop system with a linear $\L_2$-gain from disturbance input 
to the penalty variable.  This line of research is referred to as nonlinear $\mathcal{H}^\infty$ optimal
control (see for example \cite{BaBe08} \cite{HeJa99} \cite{vdSc00}).

In contrast to the explicit quantitative $\L_2$-gain design goal above, ISS was formulated as a 
qualitative robust stability property explicitly for nonlinear systems.  While there has been recent work
on computing ISS gains \cite{HJND05-TAC}, the feedback design techniques to achieve 
ISS typically rely on Lyapunov-based techniques such as control-Lyapunov functions \cite{KKK95}, \cite{SJK97}.
Consequently, while the design techniques for ISS are more easily applied to nonlinear systems,
they generally lack the pre-specified gain limits of nonlinear $\mathcal{H}^\infty$ control.

%

Sontag \cite{Sont98-SCL} investigated integral variants of the ISS property and demonstrated that
ISS is equivalent to an integral-to-integral ISS-type estimate (stated here as \eqref{eq:ISS}).  He termed
this an ``$\L_2$ to $\L_2$ property'' as, by a particular choice of the scaling functions involved,
one exactly recovers the standard definition of linear $\L_2$-gain.  In addition, Sontag observed that
by taking an integral of the input, but not the state, one obtains
a fundamentally different stability property, which he termed integral ISS (usually abbreviated to iISS).  In \cite{Sont98-SCL}
iISS is referred to as an ``$\L_2$ to $\L_\infty$ property'' and is shown to be strictly weaker
than ISS; i.e., all ISS systems are iISS but there exist iISS systems that are not ISS.  

Inspired by the nonlinear gains used in ISS-type estimates, we explicitly considered the notion of nonlinear
$\L_2$-gain \cite{DoKe-CDC08}, where the energy of the state or output penalty variable is bounded from above by a
nonlinear scaling of the energy of the input.  This generalization of linear $\L_2$-gain is intuitively
appealing as one would not {\em a priori} expect a linear bound for nonlinear systems.
In principle, the nonlinear $\L_2$-gain property applies to a wider class of systems than
does the linear $\L_2$-gain property.  Furthermore, when dealing with quantitative results,
nonlinear gains allow for tighter gain bounds, allowing for more precise results.
We subsequently developed verification \cite{DZK12-SCL}
and synthesis \cite{ZDK-IFAC11} tools for the nonlinear $\L_2$-gain property.
These tools can be seen as an extension of nonlinear $\mathcal{H}^\infty$ control.

The ISS and $\L_2$-gain approaches developed in parallel and largely independent of each
other.  A rare point of contact is the work of Gr\"{u}ne, Sontag, and 
Wirth \cite{GSW99-SCL} where, for systems of dimension different from 4 or 5, a certain qualitative equivalence was demonstrated
between global asymptotic stability of the origin and global exponential stability of
the origin (and hence $\L_2$-stability of the associated system).  Additionally, a similar
qualitative equivalence between ISS and linear $\L_2$-gain was demonstrated.  To be precise,
Gr\"{u}ne, {\it et al.\ } showed that a system with linear $\L_2$-gain is always ISS and that
given an ISS system it is possible to find a nonlinear change of coordinates so that, in the
new coordinates, the system has linear $\L_2$-gain \cite[Theorems 3, 4]{GSW99-SCL}.  In this context,  
nonlinear $\mathcal{H}^\infty$ control can be seen as a method to design ISS systems
with a prescribed ISS gain.  Design tools for attaining pre-specified ISS gains for 
nonlinear discrete-time systems were presented in \cite{HJND05-TAC}.
In Section~\ref{sec:QualEqs} we present a result similar to that of \cite[Theorems 3, 4]{GSW99-SCL}
demonstrating a qualitative equivalence between ISS and linear $\L_2$-gain (we will
make the differences precise in Theorem~\ref{thm:ISSandLinear}.)  We then present
a qualitative equivalence between iISS and nonlinear $\L_2$-gain (Theorem~\ref{thm:iISSandNonlinear}).  As a consequence,
 the synthesis and verification results for nonlinear $\L_2$-gain from \cite{ZDK-IFAC11} 
 and \cite{DZK12-SCL} can be seen as design tools for iISS systems.  To date, design tools
 for iISS systems are largely unavailable.

A natural approach to analyzing large-scale dynamical systems involves
 separating the large-scale system into several smaller
interconnected subsystems, analyzing the subsystems, and then investigating 
overall system behavior on the basis of subsystem behavior and their
interconnections.  Both ISS and $\L_2$-gain have been widely used in this context.
Consequently, it is of interest to discuss how both cascade and feedback interconnections
behave for the different stability properties.
It is immediately obvious that the cascade connection of two systems with linear $\L_2$-gain results
in an overall system with linear $\L_2$-gain.  A small-gain condition
\cite[Theorem 1]{Zame66-TAC} is sufficient to guarantee that the feedback interconnection of systems with
linear $\L_2$-gain results in an overall system with linear $\L_2$-gain.
Similarly, the cascade connection of two ISS systems is  ISS \cite{Sont89-TAC} and a small-gain condition is sufficient
to guarantee that the feedback interconnection of
two ISS systems is also ISS \cite{JTP94-MCSS}.

As we show in Section~\ref{sec:NonlinearL2gain}, when considering system interconnections,
the nonlinear $\L_2$-gain property shares many similarities with ISS and linear $\L_2$-gain.  In Section~\ref{sec:L2_cascade}
we show that the cascade of two systems with nonlinear $\L_2$-gain also has nonlinear $\L_2$-gain.
In Section~\ref{sec:L2_feedback} we show that if a small-gain condition is satisfied then the feedback interconnection
of two systems with nonlinear $\L_2$-gain also has nonlinear $\L_2$-gain.
By contrast, it is known that a cascade interconnection of iISS systems is not necessarily iISS
\cite{AAS02-SICON} and, even if a small gain condition is satisfied, a feedback interconnection of
iISS systems is not necessarily iISS.
Consequently, the aforementioned qualitative equivalence between nonlinear $\L_2$-gain and iISS
(Theorem~\ref{thm:iISSandNonlinear}) and the results of section~\ref{sec:L2_feedback} appear to contradict
known results on interconnections of iISS systems.

This apparent contradiction is resolved in Section~\ref{sec:iISS_connections} by recognizing that, while nonlinear $\L_2$-gain and iISS are
qualitatively equivalent, the relationship is asymmetric in the sense that {\em all} systems with the nonlinear
$\L_2$-gain property are iISS, while there {\em exists} a coordinate transformation for a given iISS system
so that, in the new coordinates, the system satisfies the nonlinear $\L_2$-gain property.  Consequently,
studying interconnections of iISS systems via the asserted nonlinear $\L_2$-gain qualitative equivalence requires
careful consideration of the state and input transformations used.  This consideration leads to
several sufficient conditions for stability of interconnected iISS systems similar to those
found in \cite{AAS02-SICON}, \cite{ChAn08-SCL}, and \cite{ItJi09-TAC}.

The paper is organized as follows.  In Section~\ref{sec:math_prelim} we present 
some necessary mathematical preliminaries including precise definitions of the
stability concepts of interest, as well as two key lemmas on nonlinear
changes of coordinates.  In Section~\ref{sec:QualEqs} we demonstrate essential
qualitative equivalences between the six different stability properties of interest.
In Section~\ref{sec:NonlinearL2gain} we present sufficient conditions for the 
interconnection (cascade or feedback) of systems with nonlinear \ltwo-gain to
also have nonlinear \ltwo-gain while
in Section~\ref{sec:iISS_connections} we present several sufficient 
conditions for \ltwo-stability, ISS, or iISS of interconnected systems by drawing on
the qualitatively equivalent properties developed in Section~\ref{sec:QualEqs}.
In Section~\ref{sec:Conc} we provide some concluding remarks.


\section{Preliminaries}
\label{sec:math_prelim}
We  consider systems described by ordinary differential equations of the form
\begin{equation}
	\label{eq:ODE}
	\tfrac{d}{dt}x(t) = f(x(t)), \quad x(0) \in \R^n,
\end{equation}
where $f:\R^n \rightarrow \R^n$ is locally Lipschitz.  We also
consider systems with inputs described by
\begin{equation}
	\label{eq:sys}
	\tfrac{d}{dt}x(t) = f(x(t),w(t)), \quad x(0) \in \R^n,
\end{equation}
where $f:\R^n \times \R^m \rightarrow \R^n$ is locally Lipschitz in its first argument, locally
uniformly in its second argument.  We take as the class of inputs, $\W^m$, those functions $w: \R_{\geq 0} \rightarrow \R^m$ that are 
measurable and locally essentially bounded.
We make the standing assumption that both \eqref{eq:ODE} and \eqref{eq:sys} are forward complete
(see \cite{AnSo99-SCL} for sufficient conditions).
  We will make use of the standard function classes
$\Kinf$ and $\KL$ (see \cite{Hahn67} or \cite{Kell13-sub}).  For a measurable, locally essentially bounded function
$y:\R_{\geq 0} \rightarrow \R^n$ we denote the squared two-norm by
$\twonorm{y} \doteq \int_0^t |y(\tau)|^2 d\tau$. 

\begin{remark}
	When considering $\L_2$-type properties, it is standard to take inputs from the space
	of locally $\L_2$ functions, denoted by \ltwoloc \
	 (i.e., those functions whose truncation to any finite time horizon
	is in $\L_2$).  However, a subspace of \ltwoloc, consisting of all measurable and locally
	essentially bounded functions is sufficient for what follows.  With respect to the two-norm, we exclusively
	use the truncated two-norm above which is finite for any fixed $t \in \R_{\geq 0}$ and any 
	$w \in \W^m$.  Furthermore, this class of inputs is commonly used in the ISS literature since,
	again for any fixed $t \in \R_{\geq 0}$, it guarantees that the integral of nonlinearly scaled versions
	of the input is finite over $[0,t]$ (see \cite{Sont98-SCL} or Lemma~\ref{lem:IntegrableChanges} below).
\end{remark}

\subsection{Stability Properties}
\label{sec:StabProp}
There are six stability properties that will be of interest in the sequel.  The first two 
properties are for systems without inputs \eqref{eq:ODE}, while the final four are for
systems with inputs \eqref{eq:sys}.

\begin{defn}
	System \eqref{eq:ODE} is $\alpha$-integrable if there exists $\alpha, \beta \in \Kinf$ so that
	\begin{equation}
		\label{eq:alpha_integrable}
		\int_0^t \alpha(|x(\tau)|)d\tau \leq \beta(|x(0)|), \quad \forall \, x(0) \in \R^n, \ t \in \R_{\geq 0}.
	\end{equation}
\end{defn}

In \cite{TPL02-MCSS} it was observed that $\alpha$-integrability is equivalent to uniform global asymptotic 
stability of the origin for \eqref{eq:ODE}.

\begin{defn}
	System \eqref{eq:ODE} is $\L_2$-stable if there exists $\beta \in \Kinf$ so that
	\begin{equation}
		\twonorm{x} \leq \beta(|x(0)|), \quad \forall \, x(0) \in \R^n, \ t \in \R_{\geq 0}.
	\end{equation}
\end{defn}

We observe that \ltwo-stability is a special case of $\alpha$-integrability
where $\alpha(s) = s^2$ for all $s \in \R_{\geq 0}$.

We now define four stability properties for systems with inputs.  
The first two are the well-known properties of Input-to-State
Stability (ISS) \cite{Sont89-TAC} and integral Input-to-State Stability (iISS) \cite{Sont98-SCL}.
\begin{defn}
	\label{def:ISS}
	System \eqref{eq:sys} is Input-to-State Stable (ISS) if there exist $\alpha, \beta, \sigma \in \Kinf$
	so that the estimate
	\begin{equation}
		\label{eq:ISS}
		\int_0^t \alpha(|x(\tau)|) d\tau \leq \max \left\{ \beta(|x(0)|),  \int_0^t \sigma(|w(\tau)|)d\tau \right\}
	\end{equation}
	holds for all $x(0) \in \R^n$, $w \in \W^m$, and $t \in \R_{\geq 0}$.
\end{defn}

\begin{defn}
	\label{def:iISS}
	System \eqref{eq:sys} is integral Input-to-State Stable (iISS) if there exist
	$\alpha, \beta, \gamma, \sigma \in \Kinf$ so that the estimate
	\begin{equation}
		\label{eq:iISS}
		\int_0^t \alpha \left( |x(\tau)| \right) d\tau \leq \max \left\{ \beta(|x(0)|), \gamma \left( \int_0^t \sigma(|w(\tau)|) d\tau \right) \right\}
	\end{equation}
	holds for all $x(0) \in \R^n$, $w \in \W^m$, and $t \in \R_{\geq 0}$.
\end{defn}

The final two stability properties are based on the $\L_2$-norm.

\begin{defn}
	System \eqref{eq:sys} has  linear $\L_2$-gain with transient and gain bound $\beta \in \Kinf$,
	 $\bar{\gamma} \in  \R_{\geq 0}$ if the estimate
	\begin{equation}
		\label{eq:L2gain}
		\twonorm{x} \leq \max \left\{ \beta(|x(0)|), \bar{\gamma}^2 \twonorm{w} \right\}
	\end{equation}
	holds for all $x(0) \in \R^n$, $w \in \W^m$, and $t \in \R_{\geq 0}$.
\end{defn}

\begin{defn}
	System \eqref{eq:sys} has  nonlinear $\L_2$-gain with transient and gain bound $\beta, \gamma
	\in  \Kinf$ if the estimate
	\begin{equation}
		\label{eq:nonlinearL2gain}
		\twonorm{x} \leq \max \left\{ \beta(|x(0)|), \gamma\left(\twonorm{w}\right) \right\}
	\end{equation}
	holds for all $x(0) \in \R^n$, $w \in \W^m$, and $t \in \R_{\geq 0}$.
\end{defn}

With the standing assumption that \eqref{eq:sys} is forward complete, the original definitions of ISS \cite{Sont89-TAC} and iISS \cite{Sont98-SCL} 
were shown to be equivalent to to Definition~\ref{def:ISS} and Definition~\ref{def:iISS}
in \cite{Sont98-SCL} and \cite{ASW00-DC}, respectively.


\subsection{Changes of Coordinates}
In \cite{Sontag-Survey}, Sontag asserted that
``notions of stability should be invariant under (nonlinear) changes of variables.''
In part, this derives from the fact that in order to apply various nonlinear control design methods the
system equations are usually required to be in a certain {\it normal form}.  When a system model as
given is not in the necessary normal form, a common technique is to search for a 
change of coordinates such that, in the new coordinates, the normal form is achieved and a stabilizing
control design can be undertaken.  However, unless invariance of the stability property is guaranteed
under changes of coordinates, such a stabilizing design in the new coordinates may fail to be stabilizing
in the original, possibly physically meaningful, coordinates.

\begin{defn}[\cite{Sontag-Survey}]
	A change of coordinates is a homeomorphism $T:\R^p \rightarrow \R^p$ that fixes the origin.  
	In other words, $T(\cdot)$ is continuous with a well-defined and continuous inverse 
	$T^{-1}:\R^p \rightarrow \R^p$ and such that $T(0) = 0$.
\end{defn}

If it is desirable to express the system differential equations \eqref{eq:sys} in new coordinates, then the change of 
coordinates $T(\cdot)$ must be differentiable, at least away from the origin.  However, 
the results and discussion in this paper relate to trajectory based properties, and so 
do not require differentiability of $T(\cdot)$.

The following fact was observed in \cite{Sontag-Survey} and is a useful tool for analyzing the effect of changes of coordinates
on stability properties.
\begin{lemma}
	\label{lem:BoundChanges}
	Given any change of coordinates $T:\R^p \rightarrow \R^p$, there exist $\underline{\alpha}, \overline{\alpha} \in \Kinf$
	such that
	\begin{equation}
		\underline{\alpha}(|\zeta|) \leq |T(\zeta)| \leq \overline{\alpha}(|\zeta|), \quad \forall \zeta \in \R^p.
	\end{equation}
\end{lemma}

{\it Proof:} Simply take
$ \underline{\alpha}(r) \doteq \min_{|x| \geq r} |T(x)|$ and $\overline{\alpha}(r) \doteq \max_{|x| \leq r} |T(x)|$. 
$\halmos$

An immediate consequence of the above lemma is that integrability is preserved under changes of
coordinates.
\begin{lemma}
	\label{lem:IntegrableChanges}
	Let $\xi: \R_{\geq 0} \rightarrow \R^n$ be measurable and locally essentially bounded,
	$T:\R^n \rightarrow \R^n$ be any change of coordinates, and take any $\alpha_1, \alpha_2 \in \Kinf$.
	For any $t \in \R_{\geq 0}$, the integrals
	\begin{equation}
		\label{eq:dag}
		\int_0^t \alpha_1(|\xi(\tau)|) d\tau
	\end{equation}
	and
	\begin{equation}
		\label{eq:ddag}
		\int_0^t \alpha_2(|T(\xi(\tau))|)d\tau
	\end{equation}
	exist and are finite.
\end{lemma}

{\it Proof:}
	From Lemma~\ref{lem:BoundChanges}, there exists $\overline{\alpha}_T \in \Kinf$ so that 
	$|T(\xi(\tau))| \leq \overline{\alpha}_T(|\xi(\tau)|)$ for all $\tau \in \R_{\geq 0}$.  Therefore,
	proving \eqref{eq:ddag} exists and is finite reduces to proving \eqref{eq:dag} exists
	and is finite.  To prove \eqref{eq:dag} exists and is finite we simply note that $\alpha_1 \in \Kinf$
	and the norm are both continuous functions on $\R_{\geq 0}$ and $\R^n$, respectively.  
	Consequently, if $\tau \mapsto \xi(\tau)$ is measurable
	and essentially bounded then $\tau \mapsto \alpha_1(|\xi(\tau)|)$ is also measurable and 
	locally essentially bounded, yielding the desired result.
$\halmos$

With Lemma~\ref{lem:BoundChanges} available, it is straightforward to see that $\alpha$-integrability, ISS, and iISS
satisfy Sontag's assertion that stability notions should be invariant under changes of coordinates.
Using ISS as an example, given any change of coordinates on the state $T:\R^n \rightarrow \R^n$, let the functions
$\underline{\alpha}_T, \overline{\alpha}_T \in \Kinf$ come from Lemma~\ref{lem:BoundChanges}, and define
$\Kinf$ functions $\tilde{\alpha} \doteq \alpha \circ \overline{\alpha}_T^{-1}$ and
$\tilde\beta \doteq \beta \circ \underline{\alpha}_T^{-1}$ where $\alpha, \beta \in \Kinf$ are from the ISS estimate \eqref{eq:ISS}.  
Furthermore,  let $S:\R^m \rightarrow \R^m$, be any change of coordinates on the input space, let
$\underline{\alpha}_S \in \Kinf$ come from Lemma~\ref{lem:BoundChanges}, and define the $\Kinf$ function
$\tilde\sigma \doteq \sigma \circ \underline{\alpha}_S^{-1}$ with $\sigma \in \Kinf$ from the ISS estimate
\eqref{eq:ISS}.
Define $v(t) \doteq S(w(t))$ for all $t \in \R_{\geq 0}$ and $\psi(t) \doteq T(x(t))$ for all 
$t \in \R_{\geq 0}$
Then the bounds from Lemma~\ref{lem:BoundChanges} and the ISS estimate \eqref{eq:ISS} yield
\begin{align}
	\lefteqn{ \int_0^t \tilde{\alpha}(|\psi(\tau)|)d\tau} &  \nn \\
	 &  =   \int_0^t \tilde{\alpha}\left( |T(x(\tau))| \right) d\tau
	 \leq  \int_0^t \tilde\alpha \left( \overline{\alpha}_T
		(|x(\tau)|)\right)d\tau \nonumber \\
	& =  \int_0^t \alpha(|x(\tau)|) d\tau 
	 \leq  \max \left\{ \beta(|x(0)|), \int_0^t \sigma(|w(\tau)|) d\tau \right\} \nonumber \\
	& \leq  \max \left\{ \beta\left( \underline{\alpha}_T^{-1}(|T(x(0))|)\right), 
		 \int_0^t \sigma \circ \underline{\alpha}_S^{-1}(|S(w(\tau))|) d\tau  \right\} \nn \\
	& =  \max \left\{ \tilde\beta(|\psi(0)|), \int_0^t \tilde\sigma(|v(\tau)|) d\tau \right\} . \label{eq:ISS_qual_eq}
\end{align}
In other words, the system in the new coordinates also satisfies an ISS estimate \eqref{eq:ISS} with 
functions $\tilde\alpha, \tilde\beta, \tilde\sigma \in \Kinf$ in place of $\alpha, \beta, \sigma \in \Kinf$.

Note that precisely the same argument as above holds for systems that are $\alpha$-integrable or 
which satisfy an iISS estimate \eqref{eq:iISS}.
Hence, $\alpha$-integrability, ISS, and iISS are invariant under changes of coordinates in the input and state variables.  However, none of $\L_2$-stability,  linear $\L_2$-gain, or  nonlinear $\L_2$-gain satisfy this property.  (See Examples~\ref{ex:GAS_L2} and \ref{ex:ISS_L2} in Sections \ref{sec:L2_alphaint} and \ref{sec:L2_ISS},
respectively.)
		
Similar to how the magnitude of coordinate transformations can be upper and lower bounded by functions
of class $\Kinf$, given any function of class $\Kinf$ there exist changes of coordinates that, in magnitude,
upper and lower bound this $\Kinf$ function. 
\begin{lemma}[{\cite[Lemma 2.11]{DKZ12-CDC}}]
	\label{lem:BoundKinf}
	Given any $\alpha \in \Kinf$ and $p \in \Z_{>0}$ there exist changes of coordinates $T_\ell, T_u : \R^p \rightarrow \R^p$
	such that
	\begin{equation}
		\label{eq:BoundKinf}
		|T_\ell(\zeta)| \leq \alpha(|\zeta|) \leq |T_u(\zeta)|, \quad \forall \zeta \in \R^p.
	\end{equation}
\end{lemma}

{\it Proof:} Construct the change of coordinates $T_u : \R^p \rightarrow \R^p$ with the $i^{\rm th}$
coordinate given by
\begin{equation}
	\label{eq:upper_change_def}
	T_{u,i}(\zeta) \doteq {\rm sgn} (\zeta_i) \alpha \left(|\zeta_i| \sqrt{p}\right),
\end{equation}
which is invertible by inspection.  
Continuity of the change of coordinates follows from the continuity of $\alpha \in \Kinf$ and the fact that
$\alpha(0) = 0$.
Then,
\begin{eqnarray}
	|T_u(\zeta)| & = & \sqrt{ \sum_i |T_{u,i}(\zeta)|^2} \ \geq \ \max_i |T_{u,i}(\zeta)| \nn \\
	 & = & \max_i \alpha(|\zeta_i| \sqrt{p}) 
	   =  \alpha\left(\max_i|\zeta_i| \sqrt{p} \right) 
	 \geq  \alpha(|\zeta|) . \nn
\end{eqnarray}

Similarly, construct the change of coordinates $T_\ell : \R^p \rightarrow \R^p$ with $i^{\rm th}$
coordinate
\begin{equation}
	\label{eq:lower_change_def}
	T_{\ell, i}(\zeta) \doteq \tfrac{1}{\sqrt{p}} {\rm sgn}(\zeta_i) \alpha(|\zeta_i|),
\end{equation}
which is also invertible by inspection and continuous.  Then 
\begin{eqnarray}
	\lefteqn{|T_\ell(\zeta)|  =  \sqrt{\sum_i |T_{\ell,i}(\zeta)|^2} \ \leq \ \sqrt{p} \max_i |T_{\ell,i}(\zeta)|} & &  \nn \\
	 & = & \sqrt{p} \max_i \tfrac{1}{\sqrt{p}} \alpha(|\zeta_i|) 
	  \leq  \alpha \left( \max_i |\zeta_i| \right) 
	 \leq  \alpha(|\zeta|). \nn
\end{eqnarray}
Therefore, the changes of coordinates \eqref{eq:upper_change_def} and \eqref{eq:lower_change_def}
satisfy \eqref{eq:BoundKinf}. $\halmos$

We observed above that we merely require a change of coordinates to be a homeomorphism 
since our interest herein is limited to trajectory-based properties.  However, since the subsequent
results largely rely on Lemmas \ref{lem:BoundChanges} and \ref{lem:BoundKinf}, it is worth noting
that requiring a change of coordinates to be a diffeomorphism away from the origin is not particularly
restrictive.  Indeed, with regards to Lemma~\ref{lem:BoundChanges}, the given change of coordinates may
well be a diffeormorphism.  With respect to Lemma~\ref{lem:BoundKinf}, as is evident from the proof,
both changes of coordinates can be chosen so that they inherit the regularity properties of the 
given $\Kinf$ function away from the origin.  Hence, if the given function
$\alpha \in \Kinf$ is smooth, then both changes of coordinates $T_{\ell}(\cdot)$ and
$T_{u}(\cdot)$ can be chosen to be smooth away from the origin.  Furthermore, we note that
any given function of class $\Kinf$ can be approximated with arbitrary precision by a class $\Kinf$
function that is smooth (e.g., \cite[Lemma 6]{Kell13-sub}).

%
%


\section{Qualitatively Equivalent Robust Stability Properties}
\label{sec:QualEqs}

\subsection{Qualitative Equivalence}
In the result of \eqref{eq:ISS_qual_eq},
while the comparison functions differ between the original ISS estimate \eqref{eq:ISS} and 
the ISS estimate for the new coordinates \eqref{eq:ISS_qual_eq}, the 
form of the inequality is clearly the same, and consequently \eqref{eq:ISS} and \eqref{eq:ISS_qual_eq}
are said to be {\it qualitatively equivalent}.  Where the comparison functions are the same, this
equivalence is said to be {\it quantitative}.

The notions of qualitative and quantitative equivalence 
can be extended to pairs of properties of different forms.  In particular, if a given property
implies a second property of a different form which, in turn, implies a third property
of the same form as the first, and if the first and third properties are qualitatively equivalent,
then we refer to all three properties as being qualitatively equivalent.

For example, bounds defined by sums and maximums are qualitatively 
equivalent.  
That these provide qualitatively
equivalent properties follows from the fact that for any $a, b \in \R_{\geq 0}$
\begin{equation}
	\label{eq:sum_bound}
	 a + b \leq \max\{ 2a, 2b\} \quad {\rm and} \quad 
 	\max\{a, b\} \leq a + b. 
\end{equation}

Applied to the ISS estimate \eqref{eq:ISS}, the above inequalities yield
\begin{eqnarray}
	\int_0^t \alpha(|x(\tau)|) d\tau & \leq & \max \left\{ \beta(|x(0)|), \int_0^t \sigma(|w(\tau)|)d\tau \right\} \label{eq:prop1} \\
	& \leq & \beta(|x(0)|) + \int_0^t \sigma(|w(\tau)|)d\tau \label{eq:prop2} \\
	& \leq & \max \left\{ 2\beta(|x(0)|), 2\int_0^t \sigma(|w(\tau)|)d\tau \right\} \label{eq:prop3}.
\end{eqnarray}
Clearly \eqref{eq:prop1} and \eqref{eq:prop3} are qualitatively equivalent and, by our extended notion of qualitative equivalence,
\eqref{eq:prop2} is qualitatively equivalent to \eqref{eq:prop1}.  We note that the nature of this equivalence is {\it qualitative} rather than {\it quantitative}
since the comparison function bounds are not the same due to the factor of 2 involved in the first relation in \eqref{eq:sum_bound}.


Similarly, the definitions of ISS (Definition~\ref{def:ISS})
and iISS (Definition~\ref{def:iISS}) are qualitatively equivalent to the original definitions proposed in the literature.
In particular, under the assumption of forward completeness of \eqref{eq:sys}, \cite[Theorem 1]{Sont98-SCL} demonstrated
that \eqref{eq:ISS} is qualitatively equivalent to the original definition of ISS in \cite{Sont89-TAC}; i.e.,
there exists $\gamma \in \Kinf$ and $\beta \in \KL$ so that
\begin{equation}
	\label{eq:ISS_original}
	|x(t)| \leq \max \left\{ \beta(|x(0)|,t), \sup_{\tau \in [0,t]} \gamma(|w(\tau)|) \right\}
\end{equation}
holds for all $x(0) \in \R^n$, $w \in \W^m$, and $t \in \R_{\geq 0}$. 
Similarly, again under the assumption of forward completeness of \eqref{eq:sys}, \cite[Theorem 1]{ASW00-DC} showed that \eqref{eq:iISS} is qualitatively equivalent
to the existence of $\alpha, \gamma \in \Kinf$ and $\beta \in \KL$ so that
\begin{equation}
	\label{eq:iISS_original}
	\alpha(|x(t)|) \leq \max \left\{\beta(|x(0)|,t), \int_0^t \gamma(|w(\tau)|)d\tau \right\} 
\end{equation}
holds for all $x(0) \in \R^n$, $w \in \W^m$, and $t \in \R_{\geq 0}$, as defined in \cite{Sont98-SCL}.

The original definitions of ISS, \eqref{eq:ISS_original}, and iISS, \eqref{eq:iISS_original},
possess some appealing intuitive properties.  For example,
 ISS involves bounds on system trajectories at a given time that depend
on a decaying transient term due to the initial condition ($\beta \in \KL$) as well as an additional term due to the worst-case input
up to the current time ($\gamma \in \Kinf$).  This desired property is more obvious in \eqref{eq:ISS_original}
than in the qualitatively equivalent definition \eqref{eq:ISS}.  Similarly, the fact that the input is treated in a fundamentally
different manner for integral ISS than it is for ISS is more obvious in the difference between \eqref{eq:ISS_original} and
\eqref{eq:iISS_original} than it is in the difference between \eqref{eq:ISS} and \eqref{eq:iISS}.  
On the other hand, the fact that iISS is a strictly weaker property than ISS is more obvious when examining
\eqref{eq:ISS} and \eqref{eq:iISS} than it is when examining \eqref{eq:ISS_original} and \eqref{eq:iISS_original}.
Indeed, all ISS systems are iISS since the identity is simply one possible choice of the function
$\gamma \in \Kinf$ of \eqref{eq:iISS}.  Furthermore, since there are many $\Kinf$ functions which are not the identity,
iISS possibly encompasses a larger class of systems.  That this is in fact the case was shown in
\cite{Sont98-SCL}.  Therefore, we see that by examining qualitatively equivalent 
ISS properties and their relationships to qualitatively equivalent iISS properties, we gain a clearer understanding
of the relationship between ISS and iISS systems.
Furthermore, as is evident in the sequel, the ISS and iISS definitions given by \eqref{eq:ISS} and \eqref{eq:iISS},
respectively, are better suited to clarifying the relationship between these properties and the  $\L_2$-gain properties \eqref{eq:L2gain} and \eqref{eq:nonlinearL2gain}.

With the notion of qualitative equivalence established, the remainder of this section is concerned with
establishing that $\alpha$-integrability, ISS, and iISS are qualitatively equivalent, via a change of
coordinates, to $\L_2$-stability,  linear $\L_2$-gain, and  nonlinear $\L_2$-gain, respectively.  It is important to
note that these equivalences are not quantitative; e.g., the iISS-gain $\gamma \in \Kinf$ of \eqref{eq:iISS} is not, in general,
the  nonlinear $\L_2$-gain $\gamma \in \Kinf$ of \eqref{eq:nonlinearL2gain}.

\subsection{$\L_2$-stability and $\alpha$-integrability}
\label{sec:L2_alphaint}
\begin{thm}
	\label{thm:alpha_l2stable}
	If system \eqref{eq:ODE} is $\L_2$-stable then it is $\alpha$-integrable.
	Conversely, if system \eqref{eq:ODE} is $\alpha$-integrable then there
	exists a change of coordinates for the state such that the system in the
	new coordinates is $\L_2$-stable.
\end{thm}

{\it Proof:} That $\L_2$-stability implies $\alpha$-integrability is obvious by inspection
since $\alpha(s) \doteq s^2$ for all $s \in \R_{\geq 0}$ is of class $\Kinf$.
In order to prove the converse, assume we have $\alpha,\beta \in \Kinf$ so that
\eqref{eq:alpha_integrable} is satisfied.  Apply Lemma~\ref{lem:BoundKinf} to
$\alpha^{\frac{1}{2}} \in \Kinf$ with $p = n$ to obtain $T :\R^n \rightarrow \R^n$ so that
$|T(x)| \leq \alpha^{\frac{1}{2}}(|x|)$, for all $x \in \R^n$. 
Lemma~\ref{lem:BoundChanges} implies the existence of $\alpha_T \in \Kinf$ such that
$\alpha_T(|x|) \leq |T(x)|$ for all $x \in \R^n$.
Defining the new coordinates $\xi \doteq T(x)$, we then see that, for all $\xi(0) \in \R^n$ and $t \in \R_{\geq 0}$,
\begin{eqnarray}
	\twonorm{\xi} & = & \int_0^t |\xi(\tau)|^2 d\tau =  \int_0^t |T(x(\tau))|^2 d\tau \nn \\
	& \leq & \int_0^t \alpha(|x(\tau)|) d\tau \ \leq \ \beta(|x(0)|) \nn \\
	 & \leq & \beta \circ \alpha_T^{-1}(|T(x(0))|) = \beta \circ \alpha_T^{-1}(|\xi(0)|) \nn
\end{eqnarray}
so that, in the new coordinates, system \eqref{eq:ODE} is $\L_2$-stable.
$\halmos$

\begin{example}
\label{ex:GAS_L2}
The origin can be shown to be globally asymptotically stable for 
\begin{equation}
	\label{eq:ex1_sys}
	\tfrac{d}{dt}x(t) = -x(t)^3 , \qquad x(0) \in \R
\end{equation}
by using the Lyapunov function $V(x) = \frac{1}{2}x^2$.
Consequently, by the observation in \cite{TPL02-MCSS}, \eqref{eq:ex1_sys} is
$\alpha$-integrable.
The solution of \eqref{eq:ex1_sys} is
\begin{equation}
	x(t) = \frac{x(0)}{\sqrt{1+2x(0)^2t}} , \quad \forall \, x(0) \in \R, t \in \R_{\geq 0} \nn
\end{equation}
so that
\begin{equation}
	\twonorm{x} = \frac{1}{2} \log \left( 1 + 2x(0)^2t \right) \nn
\end{equation}
and hence there is no $\beta \in \Kinf$ such that
$\twonorm{x} \leq \beta(|x(0)|)$. 
In other words, while \eqref{eq:ex1_sys} is $\alpha$-integrable, it is not $\L_2$-stable.
However, Theorem~\ref{thm:alpha_l2stable} states that there exists a change of coordinates
so that, in the new coordinates, the system is $\L_2$-stable.  Let
\begin{equation}
	\label{eq:ex_change}
	z = T(x) \doteq x\exp\left(-\frac{1}{2x^2}\right), \quad x \in \R \backslash \{0\} ,
\end{equation}
and $T(0) = 0$.  This change of coordinates is a homeomorphism on $\R$ and
a diffeomorphism on $\R \backslash \{0\}$.
It is straightforward to write the system equation in the new coordinates as
\begin{equation}
	\label{eq:z_L2}
	\tfrac{d}{dt}z(t) = -z(t) \left(1 + |T^{-1}(z(t))|^2\right), \quad z(0) \in \R \backslash \{0\} 
\end{equation}
from which it follows that $|z(t)| \leq |z(0)|\exp(-t)$ for all $z(0) \in \R$ and $t \in \R_{\geq 0}$.
Consequently, $\twonorm{z(t)} \leq \frac{1}{2} |z(0)|^2$, and hence the system in the new
coordinates is $\L_2$-stable.

We also observe that this example demonstrates that $\L_2$-stability is not invariant under
changes of coordinates.  This follows from the fact that $T^{-1}(\cdot)$ is a change of coordinates
that transforms the $\L_2$-stable system \eqref{eq:z_L2} to the system \eqref{eq:ex1_sys} that is not
$\L_2$-stable.  
\end{example}

\subsection{Linear $\L_2$-gain and ISS}
\label{sec:L2_ISS}
The following theorem is similar to \cite[Theorem 4]{GSW99-SCL} demonstrating a qualitative
equivalence between ISS and (linear) \ltwo-gain, with a few key differences.

\begin{thm}
	\label{thm:ISSandLinear}
	If system \eqref{eq:sys} has  linear $\L_2$-gain then system \eqref{eq:sys} is ISS.  Conversely,
	for any $\bar\gamma^2 \in \R_{>0}$,
	if system \eqref{eq:sys} is ISS then there exist changes of coordinates for both the input
	and state such that the system in the new coordinates has  linear $\L_2$-gain $\bar{\gamma}^2$.
\end{thm}

In \cite[Theorem 4]{GSW99-SCL}, a change of coordinates is constructed such that the
bounding term related to the initial condition in \eqref{eq:L2gain} can be taken as $\beta = {\rm Id}$.
Obtaining this result relies on the level sets of an appropriate Lyapunov function being homeomorphic
(or diffeomorphic) to spheres and, as a consequence, \cite[Theorem 4]{GSW99-SCL} does not hold for dimensions
$n=4,5$.  Theorem~\ref{thm:ISSandLinear} above has no such restriction at the expense of not being able to
choose {\it a priori} the function $\beta \in \Kinf$.

We observe that one can arbitrarily set the gain parameter by
appropriate choice of the change of coordinates on the input variable.
The ability to fix the \ltwo-gain parameter in Theorem~\ref{thm:ISSandLinear} is analogous to
the ability to fix the decay rate in Sontag's lemma on $\KL$-estimates \cite[Proposition 7]{Sont98-SCL}
(also \cite[Lemma 7]{Kell13-sub}).  That is, for a given function $\beta \in \KL$ and a desired
decrease rate $\lambda \in \R_{>0}$, there exist $\alpha_1,\alpha_2 \in \Kinf$ such that
$\alpha_1(\beta(s,t)) \leq \alpha_2(s)\exp(-\lambda t)$ for all $s,t \in \R_{\geq 0}$. 

{\it Proof of Theorem~\ref{thm:ISSandLinear}:}
The proof of the first statement in Theorem~\ref{thm:ISSandLinear}
is straightforward since  linear $\L_2$-gain \eqref{eq:L2gain} is a 
special case of
the ISS estimate \eqref{eq:ISS} where the comparison functions in the ISS definition are simply
$\alpha(s) =  s^2$ and  $\sigma(s) = \bar{\gamma}^2 s^2$, for all $s \in \R_{\geq 0}$.

To show the converse statement of Theorem~\ref{thm:ISSandLinear}, 
suppose that system \eqref{eq:sys} is ISS so that \eqref{eq:ISS} is satisfied
with functions $\alpha, \beta, \gamma \in \Kinf$.  For the function $\alpha^{\frac{1}{2}} \in \Kinf$, with $p=n$ Lemma~\ref{lem:BoundKinf}
yields the existence of a change of coordinates $T:\R^n \rightarrow \R^n$ so that
\begin{equation}
	\label{eq:lower_alpha}
 	|T(x)| \leq \alpha^{\frac{1}{2}}(|x|), \quad \forall x \in \R^n .
\end{equation}
For the change of coordinates $T(\cdot)$, let $\underline{\alpha} \in \Kinf$ come from Lemma~\ref{lem:BoundChanges}
so that
\begin{equation}
	\label{eq:under_alpha}
	\underline{\alpha}(|x|) \leq |T(x)|, \quad \forall x \in \R^n .
\end{equation}
For the function $\bar\gamma^{-1} \sigma^{\frac{1}{2}} \in \Kinf$, with $p = m$ Lemma~\ref{lem:BoundKinf}  yields 
the existence of a change of coordinates
$S:\R^m \rightarrow \R^m$ such that
\begin{equation}
	\label{eq:upper_sigma}
 	\bar\gamma^{-1}\sigma^{\frac{1}{2}}(|w|) \leq |S(w)|, \quad \forall w \in \R^m. 
\end{equation}

Combining \eqref{eq:lower_alpha}, \eqref{eq:ISS}, \eqref{eq:under_alpha}, and  \eqref{eq:upper_sigma} we have,
for all $x(0) \in \R^n$, $w \in \W^m$, and $t \in \R_{\geq 0}$,
\begin{eqnarray}
	\lefteqn{ \twonorm{T(x)}  =  \int_0^t |T(x(\tau))|^2 d\tau 
	 \leq  \int_0^t \alpha(|x(\tau)|)d\tau} & & \nn \\
	 & \leq & \max \left\{ \beta(|x(0)|), \bar\gamma^2 \int_0^t \bar\gamma^{-2}\sigma(|w(\tau)|) d\tau \right\} \nonumber \\
	& \leq & \max \left\{ \beta\left( \underline{\alpha}^{-1}(|T(x(0))|) \right), \bar\gamma^2 \int_0^t |S(w(\tau))|^2 d\tau \right\} \nonumber \\
	& = & \max \left\{ \tilde\beta(|T(x(0))|), \bar\gamma^2\twonorm{S(w)} \right\} . \nn 
\end{eqnarray}
where $\tilde\beta := \beta \circ \underline{\alpha}^{-1} \in \Kinf$.  Since the input
$w(\cdot)$ is measurable and locally essentially bounded, Lemma~\ref{lem:IntegrableChanges} 
yields that $|S(w(\cdot))|^2$ is also measurable and locally essentially bounded, and hence the 
final two input-dependent terms above are well-defined.
In other words, in the state coordinates defined by $T(\cdot)$ and the input coordinates defined by
$S(\cdot)$, the system has  linear $\L_2$-gain with transient and gain bound $\tilde\beta \in \Kinf$,
$\bar\gamma^2 \in \R_{>0}$.
$\halmos$

\begin{example}
	\label{ex:ISS_L2}
	Consider the system \eqref{eq:ex1_sys} of the previous example augmented with
	an input; i.e.,
	\begin{equation}
		\label{eq:ISS_ex}
		\tfrac{d}{dt} x(t) = -x(t)^3 + w(t), \quad x(0) \in \R, \ w \in \W^1.
	\end{equation}
	Define $V(x) = \frac{1}{2}x^2$ for all $x \in \R$ and observe that
	$|x| > |w|^{1/3}$ implies $\tfrac{d}{dt}V(x(t)) < 0$.
	Therefore, $V(\cdot)$ is an ISS-Lyapunov function and, consequently, \eqref{eq:ISS_ex}
	is ISS \cite[Theorem 1]{SoWa95-SCL}.  However, by setting $w \equiv 0$, we can repeat the argument
	of Example~\ref{ex:GAS_L2} to see that \eqref{eq:ISS_ex} cannot have linear (or in fact
	nonlinear) $\L_2$-gain.
	However, as indicated by Theorem~\ref{thm:ISSandLinear}, there exists a change of
	coordinates so that, in the new coordinates, the system \eqref{eq:ISS_ex} has linear $\L_2$-gain.  Using
	the same change of coordinates \eqref{eq:ex_change} as in Example~\ref{ex:GAS_L2}, 
	we see that \eqref{eq:ISS_ex} becomes
	\begin{eqnarray}
		\tfrac{d}{dt} z(t) & = & -z(t) \left(1 + T^{-1}(z(t))^2 \right) \nn \\
			& & \ + \! \exp\left(- \frac{1}{2x(t)^2}\right) \left(1 + \frac{1}{x(t)^2} \right) w(t). 	\label{eq:z_ISS}
	\end{eqnarray}
	The inequality $1 - \frac{1}{s} \leq \log s$, $s \in \R_{>0}$, implies that
	$\exp\left(-\frac{1}{2x^2}\right) \leq \frac{2x^2}{2x^2 + 1}$
	and hence the term multiplying the input is bounded from above by $2$.  Consequently, 
	\begin{eqnarray*} 
		\tfrac{d}{dt}|z(t)| & \leq & -|z(t)| \left(1 + |T^{-1}(z(t))|^2 \right) + 2 |w(t)| \\
		& \leq & -|z(t)| + 2|w(t)| 
	\end{eqnarray*}
	and hence the system in the new coordinates has a linear $\L_2$-gain of 2.
	
	As in Example~\ref{ex:GAS_L2}, the change of coordinates $T^{-1}(\cdot)$ takes a system with
	linear $\L_2$-gain to one that is ISS but which has neither linear nor nonlinear $\L_2$-gain. 
	In this regard, with respect to Sontag's assertion that stability notions should be
	invariant under nonlinear changes of coordinates  \cite{Sontag-Survey}, linear $\L_2$-gain is not a ``good"
	notion of robust stability.  In this case, maintaining robust stability is not an issue but achieving robust
	performance may be.  In particular, if a feedback design is performed in transformed coordinates in
	order to achieve a particular linear $\L_2$-gain, there is no guarantee that system in the original,
	probably physically meaningful, coordinates will satisfy {\it any} linear $\L_2$-gain.  This is not to say
	that $\L_2$-gain is somehow an inappropriate design goal in general, as the literature demonstrates
	it has been highly successful, but that care must be taken when $\L_2$-gain techniques are coupled
	with the use of coordinate transformations.
	
	Finally, we note that Theorem~\ref{thm:ISSandLinear} suggests a method for designing Input-to-State
	Stabilizing controllers based on finding a change of coordinates so that, in the new coordinates, one can
	construct a feedback stabilizer achieving a linear $\L_2$-gain.  In the original coordinates, this then 
	provides a feedback stabilizer rendering the system ISS.  
\end{example}

\subsection{Nonlinear $\L_2$-gain and iISS}
The relationship between iISS and nonlinear $\L_2$-gain is similar to that between ISS and linear
$\L_2$-gain. 
\begin{thm}
	\label{thm:iISSandNonlinear}
	If system \eqref{eq:sys} has  nonlinear $\L_2$-gain then system \eqref{eq:sys} is iISS.
	Conversely, if system \eqref{eq:sys} is iISS then there exist changes of coordinates for
	both the input and state such
	that the system expressed in the new coordinates has  nonlinear $\L_2$-gain.
\end{thm}

One critical difference between Theorem~\ref{thm:ISSandLinear} and 
Theorem~\ref{thm:iISSandNonlinear} is in the converse statement where, in Theorem~\ref{thm:ISSandLinear},
one can choose the linear \ltwo-gain, $\bar\gamma \in \R_{>0}$, arbitrarily.  By contrast, it is
not possible to set the \ltwo-gain function in the converse statement of Theorem~\ref{thm:iISSandNonlinear}.
However, we can introduce a scaling factor inside the gain function as follows:
\begin{prop}
	\label{prop:scale_nonlinear}
	Fix $\lambda \in \R_{>0}$.  If the system \eqref{eq:sys} is iISS then there exist changes of coordinates
	$S:\R^m \rightarrow \R^m$ and $T:\R^n \rightarrow \R^n$ for the input and state, respectively,
	such that the system in the new coordinates satisfies
	\begin{equation}
		\twonorm{T(x)} \leq \max \left\{ \tilde{\beta}(|T(x(0))|), \gamma\left(\lambda \twonorm{S(w)}\right) \right\}
	\end{equation}
	for all $x(0) \in \R^n$, $w \in \W^m$, and $t \in \R_{\geq 0}$.
\end{prop}

{\it Proof of Theorem~\ref{thm:iISSandNonlinear} and Proposition~\ref{prop:scale_nonlinear}:}
As with Theorem~\ref{thm:ISSandLinear}, the proof of the first statement in Theorem~\ref{thm:iISSandNonlinear}
is straightforward since
the  nonlinear $\L_2$-gain estimate \eqref{eq:nonlinearL2gain} is an iISS estimate
\eqref{eq:iISS} with the functions $\alpha, \sigma \in \Kinf$ given by
$\alpha(s) = \sigma(s) = s^2$ for all $s \in \R_{\geq 0}$.

The proof of the converse statement of Theorem~\ref{thm:iISSandNonlinear} is a special case 
of the proof of Proposition~\ref{prop:scale_nonlinear} with $\lambda = 1$ and follows the same
argument as above for the converse statement of Theorem~\ref{thm:ISSandLinear}.  
With the function $\alpha \in \Kinf$ from \eqref{eq:iISS} we again use the state change of coordinates
 \eqref{eq:lower_alpha}  and the bound \eqref{eq:under_alpha}.
From Lemma~\ref{lem:BoundKinf}, with $p = m$,  we obtain a change of coordinates for the input satisfying
$\lambda^{-\frac{1}{2}}\sigma^{\frac{1}{2}}(|w|) \leq |S(w)|$, for all $w \in \R^m$.
We then obtain a  nonlinear $\L_2$-gain
estimate as follows:
\begin{align}
	\lefteqn{\twonorm{T(x)}  =  \int_0^t |T(x(\tau))|^2 d\tau
	 \leq  \int_0^t \alpha(|x(\tau)|)d\tau } &  \nonumber \\
	& \leq  \max \left\{ \beta(|x(0)|), \gamma \left( \lambda \int_0^t \lambda^{-1} \sigma(|w(\tau)|) d\tau \right) \right\} \nonumber \\
	& \leq  \max \left\{ \beta\left( \underline{\alpha}^{-1}(|T(x(0))|) \right), \gamma \left( \lambda \int_0^t |S(w(\tau))|^2 d\tau \right) \right\} \nonumber \\
	& =  \max \left\{ \tilde\beta(|x(0)|), \gamma \left( \lambda \twonorm{S(w)} \right) \right\} \nonumber ,
\end{align}
for all $x(0) \in \R^n$, $w \in \W^m$, and $t \in \R_{\geq 0}$. 
$\halmos$

\begin{example}
	\label{ex:iISS_nonlinearL2}
	Consider the scalar bilinear system 
	\begin{equation}
		\tfrac{d}{dt} x(t)
		 = - x(t) + x(t)\, w(t)\,, \quad x(0) \in \R, \ w \in \W^1 \,.
	\label{eq:bilinear}
\end{equation}
The iISS-Lyapunov function $V(x) \doteq \log(1+x^2)$ can be used to show that
\eqref{eq:bilinear} is iISS \cite[Theorem 1]{ASW00-TAC}.  
That \eqref{eq:bilinear} is not ISS can be seen by taking the
constant input $w(t) = 2$ for all $t \in \R_{\geq 0}$.  We now proceed to demonstrate that \eqref{eq:bilinear}
satisfies the nonlinear $\L_2$-gain property \eqref{eq:nonlinearL2gain} but not
the linear $\L_2$-gain property \eqref{eq:L2gain}.  As a consequence, just as there are iISS systems
which are not ISS, there are systems with nonlinear $\L_2$-gain which do not admit a linear
$\L_2$-gain.

Setting $Q(x) \doteq \tfrac{1}{2} \, x^2$, \eqref{eq:bilinear} implies that
\begin{align}
	\frac{1}{Q(x(\sigma))}
	\frac{dQ(x(\sigma))}{d\sigma}
	& = -2 + 2\, w(\sigma)
	\label{eq:ODE-Q}
\end{align}
for all $\sigma\in[0,s]$, $s\in\R_{\ge 0}$. Fix any $x(0)\in\R$, $t\in\R_{\ge 0}$, and any $w\in\W^1$. Integrating 
\eqref{eq:ODE-Q} over $[0,s]$, $s\in[0,t]$, yields
\begin{align}
	\log \left( \frac{Q(x(s))}{Q(x(0))} \right) 
	& = -2\, s + 2 \int_0^s w(\sigma) \, d\sigma
	\nn
\end{align}
or
	$|x(s)|^2
	 = 
	|x(0)|^2 \, \exp \left( -2\, s +  2 \int_0^s w(\sigma) \, d\sigma \right)$.
Hence,
\begin{align}
	\lefteqn{\|x\|_{\Ltwo[0,t]}^2
	 \le |x(0)|^2 \int_0^t \exp\left( -2\, s + \int_0^s 2 \, |w(\sigma)| \, d\sigma \right) ds} &
	\nn\\
	& \le |x(0)|^2 \int_0^t \exp\left( -s + \int_0^t |w(\sigma)|^2 \, d\sigma \right) ds \nn\\
	& = |x(0)|^2 \left( \int_0^t \exp(-s) \, ds \right) \exp\left(\|w\|_{\Ltwo[0,t]}^2 \right)
	\nn\\
	& = |x(0)|^2 \left( 1 - \exp(-t) \right)  \exp\left(\|w\|_{\Ltwo[0,t]}^2 \right) \nn\\
	& \le |x(0)|^2 \, \exp\left(\|w\|_{\Ltwo[0,t]}^2 \right)
	\nn\\
	& = |x(0)|^2 + |x(0)|^2 \left( \exp\left(\|w\|_{\Ltwo[0,t]}^2 \right) - 1 \right) \nn \\
	& \le |x(0)|^2 + \tfrac{1}{2}\, |x(0)|^4 + \tfrac{1}{2} \left( \exp\left(\|w\|_{\Ltwo[0,t]}^2 \right) - 1 \right)^2 .
	\label{eq:nonlinear-L2-bound}
\end{align}
Define the comparison functions
\begin{align}
	\beta(s)
	& \doteq s^2 + \tfrac{1}{2}\, s^4 \quad {\rm and} \quad
	\gamma(s) \doteq \tfrac{1}{2} \left( \exp(s) - 1 \right)^2
	\label{eq:hat-bounds}
\end{align}
for all $s\in\R_{\ge 0}$. Note that $\beta,\, \gamma\in\Kinf$. As $x(0)\in\R$, $t\in\R_{\ge 0}$, and $w\in\W^1$ in \eqref{eq:nonlinear-L2-bound} are all arbitrary, this implies that the nonlinear {\ltwo}-gain property \eqref{eq:nonlinearL2gain} holds with transient and gain bound $\beta, \gamma \in \Kinf$.

It may also be shown that system \eqref{eq:bilinear} cannot satisfy the linear {\ltwo}-gain property \eqref{eq:L2gain}. To this end, suppose that \eqref{eq:L2gain} holds with some transient and gain bound $\hat\beta \in \Kinf$,
$\bar\gamma \in \R_{\ge 0}$. Select the initial state $x(0)\in\R$ such that $|x(0)| = \hat\beta^{-1}(1) \ne 0$, and fix any  $t^{*} \in\R_{\geq 0}$ sufficiently large such that
\begin{align}
	1 + 2\, \bar\gamma^2 \, t^*
	& < \tfrac{1}{4} \left( \hat\beta^{-1}(1) \right)^2 \left( \exp(2\, t^*) -1 \right)\,.
	\label{eq:time-horizon}
\end{align}
(Note that such a $t^* \in\R_{\ge 0}$ always exists.) Select the input $w(s)  = \overline{w}(s) = 2$ for all $s\in[0,t^*]$. 
By inspection of \eqref{eq:bilinear}, $x(s) = x(0)\, \exp(s)$ for all $s\in[0,t^*]$. Hence, 
\begin{align}
	\|x\|_{\Ltwo[0,t^*]}^2
	& = \tfrac{1}{2} |x(0)|^2 \left( \exp(2\, t^*) - 1 \right) \nn \\
	& > \tfrac{1}{4} \left( \hat\beta^{-1}(1) \right)^2 \left( \exp(2\, t^*) -1 \right)
	\nn\\
	& > 1 + 2\, \bar\gamma^2 \, t^*
	 = \hat\beta(|x(0)|) + \bar\gamma^2 \|\overline{w}\|_{\Ltwo[0,t^*]}^2\,,
	\nn
\end{align}
where the second inequality above is as per \eqref{eq:time-horizon}. That is, there exist $x(0)\in\R$, $t=t^*\in\R_{\ge 0}$, and $w\in\W^1$ such that the linear {\ltwo}-gain property \eqref{eq:L2gain} with transient and gain bound  $\hat\beta \in \Kinf$, $\bar\gamma \in \R_{\ge 0}$ is violated. Furthermore, as $\hat\beta \in \Kinf$, $\bar\gamma \in\R_{\ge 0}$ are arbitrary, it follows immediately that the linear {\ltwo}-gain property \eqref{eq:L2gain} can never hold for system \eqref{eq:bilinear}.
\end{example}

\begin{remark}
We recall that \cite[Proposition 6]{Sont98-SCL} demonstrated that if system \eqref{eq:sys} satisfies
the iISS estimate \eqref{eq:iISS_original} and if the input satisfies $\int_0^t \gamma(|w(\tau)|) d\tau < \infty$
for all $t \in \R_{\geq 0}$, 
then system trajectories satisfy $x(t) \rightarrow 0$ as $t \rightarrow \infty$ 
for all $x \in \R^n$.  The obvious analogue
of this condition for the iISS estimate \eqref{eq:iISS} requires that $\int_0^t \sigma(|w(\tau)|)d\tau < \infty$ for all $t \in \R_{\geq 0}$
as this then implies 
\[ \gamma\left(\int_0^t \sigma(|w(\tau)|)d\tau \right) < \infty \]
and consequently 
\[ \int_0^t \alpha(|x(\tau)|)d\tau < \infty. \]
Finally, since $\alpha \in \Kinf$, we see that, for all $x(0) \in \R^n$, $x(t) \rightarrow 0$ 
as $t \rightarrow \infty$.
With this fact and Theorem~\ref{thm:iISSandNonlinear} we immediately see that if  system \eqref{eq:sys}
has the nonlinear $\L_2$-gain property \eqref{eq:nonlinearL2gain} then, for all $x(0) \in \R^n$, system trajectories satisfy $x(t) \rightarrow 0$
as $t \rightarrow \infty$. 
\end{remark}

\begin{remark}
The notion of $\L_2$-gain is usually stated as an {\it input-output} stability property.  If \eqref{eq:sys}
is augmented with a continuous output mapping $h:\R^n \rightarrow \R^p$ for some $p \in \Z_{>0}$
such that there exist $\underline{\alpha}_h, \overline{\alpha}_h \in \Kinf$ satisfying
	\begin{equation}
		\label{eq:ProperOutput}
		\underline{\alpha}_h(|\xi|) \leq |h(\xi)| \leq \overline{\alpha}_h(|\xi|), \quad \forall \xi \in \R^n \nn
	\end{equation}
then the previous equivalences in Theorems~\ref{thm:ISSandLinear} and \ref{thm:iISSandNonlinear}
can be shown to hold in an input-output sense.

It is unknown if the ISS equivalences, \eqref{eq:ISS} and \eqref{eq:ISS_original},
and iISS equivalences, \eqref{eq:iISS} and \eqref{eq:iISS_original}, still hold in the case of outputs that do
not satisfy \eqref{eq:ProperOutput}.  However, in \cite{KWD13-NOLCOS} it
was shown that dissipative-form and implication-form ISS-Lyapunov functions are not
equivalent in the absence of \eqref{eq:ProperOutput}.  As this equivalence
is used in the proof of \cite[Theorem 1]{Sont98-SCL}, it is possible that in the input-output case Definitions~\ref{def:ISS}
and \ref{def:iISS} are not qualitatively equivalent to the original definitions of ISS \eqref{eq:ISS_original}
and iISS \eqref{eq:iISS_original}, respectively.
Consequently, in the absence of \eqref{eq:ProperOutput} the results of Theorems~\ref{thm:ISSandLinear} and \ref{thm:iISSandNonlinear}
may not generalize to the input-output case.
\end{remark}

\section{Interconnections of Systems with Nonlinear \ltwo-gain}
\label{sec:NonlinearL2gain}
  In this section we will show that the cascade interconnection
of two systems with nonlinear \ltwo-gain itself has nonlinear \ltwo-gain.  We also present
a small-gain theorem for the feedback interconnection of systems with nonlinear \ltwo-gain
that guarantees  nonlinear \ltwo-gain for the interconnected system.  Later, in 
Section~\ref{sec:iISS_connections} we relate these results to those known to hold for iISS systems.
Here, we specifically consider two systems
\begin{eqnarray}
	\Sigma_1 & : & \tfrac{d}{dt}x_1(t) = f_1(x_1(t),w_1(t)) \label{eq:sys1} \\
	\Sigma_2 & : & \tfrac{d}{dt}x_2(t) = f_2(x_2(t),w_2(t)) \label{eq:sys2}
\end{eqnarray}
where $x_i(0) \in \R^{n_i}$, $w_i \in \W^{m_i}$, $i=1,2$, and each satisfying
\begin{align}
	\twonorm{x_1} & \leq  \max \left\{ \beta_1(|x_1(0)|), \gamma_1\left(\twonorm{w_1}\right) \right\}, \label{eq:z1_L2bnd} \\
	\twonorm{x_2} & \leq  \max \left\{ \beta_2(|x_2(0)|), \gamma_2\left(\twonorm{w_2}\right) \right\}, \label{eq:z2_L2bnd}
\end{align}
for all $x_i(0) \in \R^{n_i}$, $w_i \in \W^{m_i}$, and $t \in \R_{\geq 0}$.

A necessary prerequisite for our results on interconnecting systems is the following weak triangle inequality
 from \cite{JTP94-MCSS} (see also \cite[Lemma 4]{Kell13-sub}):
\begin{lemma}
	\label{lem:weak_tri_ineq}
	For any  $\gamma \in \Kinf$, any  $\rho \in \Kinf$ such that $\rho - {\rm Id} \in \Kinf$,
	and $a, b \in \R_{\geq 0}$,
	\begin{equation}
		\label{eq:max_tri_ineq}
		\gamma (a+b) \leq \max \left\{\gamma \circ \rho(a), \gamma \left( \rho \circ (\rho-{\rm Id})^{-1} (b) \right) \right\}.
	\end{equation}
\end{lemma}

We note that the above inequality is a generalization of the 
weak triangle inequality in \cite{Sont89-TAC}; i.e., for any $\gamma \in \K$, $a, b \in \Rnn$,
\begin{eqnarray}
	\gamma(a + b) & \leq & \max \{ \gamma(2a), \gamma(2b) \}.
	\label{eq:simple_tri_ineq} 
\end{eqnarray}

\begin{lemma}
	\label{clm:LowerBoundKSum}
	Given $\alpha_1,\alpha_2 \in \Kinf$, there exists $\alpha \in \Kinf$ so that,
	for all $s_1,s_2 \in \R_{\geq 0}$, $\alpha(s_1 + s_2) \leq \alpha_1(s_1) + \alpha_2(s_2)$.
\end{lemma}

{\it Proof:} Define 
 $\displaystyle \alpha(s) \doteq \min \left\{\alpha_1\left(\tfrac{1}{2} s \right), \ \alpha_2\left(\tfrac{1}{2} s\right)\right\}$ 
for all $s \in \R_{\geq 0}$.
Then,
\begin{multline*}
	\alpha(s_1 + s_2)  \leq  \alpha(2s_1) + \alpha(2s_2) \\
	  =  \min\left\{\alpha_1(s_1), \ \alpha_2(s_1) \right\} +
		\min\left\{\alpha_1(s_2), \ \alpha_2(s_2) \right\} \\
	 \leq  \alpha_1(s_1) + \alpha_2(s_2) . \hspace*{1.0in} \blacksquare
\end{multline*}

The following is a consequence of the definition of the $\L_2$-norm, the triangle inequality, and Young's inequality.
\begin{lemma}
	\label{clm:Young}
	For any $\varepsilon>0$ and for all $a, b \in \W^m$,
	\begin{equation}
		\twonorm{a+b} \leq (1+\varepsilon^2)\twonorm{a} 
		+ \left( 1 + \tfrac{1}{\varepsilon^2} \right) \twonorm{b}. \nn
	\end{equation}
\end{lemma}


\subsection{Cascade Interconnection}
\label{sec:L2_cascade}
\begin{figure}[h]
\centering
\includegraphics[width=6.0cm]{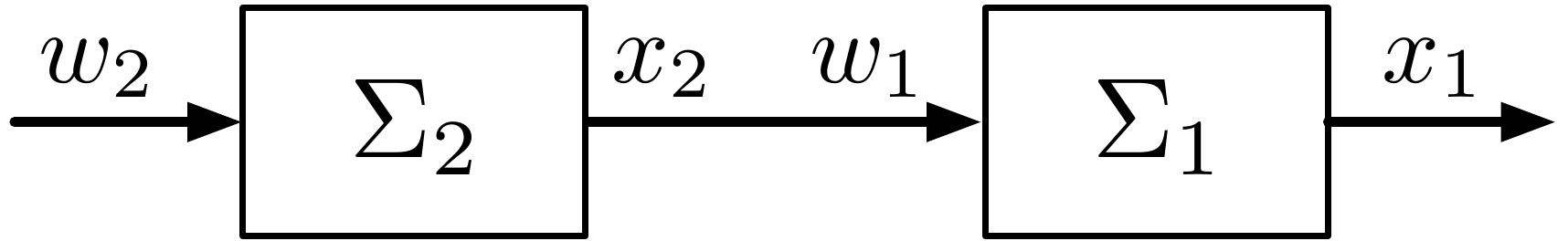}
\caption{Cascade Interconnection}
\label{fig:CascadeDiag}
\end{figure}

We first examine the cascade interconnection of \eqref{eq:sys1} and \eqref{eq:sys2}
with the interconnection $w_1 = x_2$ (requiring $m_1 = n_2$) as shown in Figure~\ref{fig:CascadeDiag}.  

\begin{prop}
	\label{prop:L2Cascade}
	Suppose that systems \eqref{eq:sys1} and \eqref{eq:sys2} satisfy the  nonlinear
	\ltwo-gain properties \eqref{eq:z1_L2bnd} and \eqref{eq:z2_L2bnd}, respectively.
	Then there exist $\beta,\gamma \in \Kinf$ such that the system given 
	by the cascade interconnection defined by $w_1 = x_2$ satisfies
	\begin{equation}
		\twonorm{x} \leq \max \left\{ \beta\left(|x(0)|\right), \gamma\left(\twonorm{w_2}\right) \right\}
	\end{equation}
	for all $x(0) \doteq [x_1^T(0) \, x_2^T(0)]^T \in \R^{n_1+n_2}$, $w_2 \in \W^{m_2}$, and $t \in \R_{\geq 0}$,
	and where $x(t) \doteq [x_1^T(t) \, x_2^T(t)]^T$ for all $t \in \R_{\geq 0}$.
\end{prop}
{\it Proof:}   Using
 bounds \eqref{eq:z1_L2bnd} and \eqref{eq:z2_L2bnd}, and the interconnection
constraint $w_1 = x_2$ we have
\begin{align}
	\lefteqn{\twonorm{x_1}  \leq  \max \left\{ \beta_1(|x_1(0)|), \gamma_1\left(\twonorm{w_1}\right) \right\} } & \nn\\
	& =  \max\left\{ \beta_1(|x_1(0)|), \gamma_1\left(\twonorm{x_2}\right) \right\}  &  \nonumber \\
	& \leq  \max\left\{ \beta_1(|x_1(0)|), \rule{0pt}{13pt}\right. \nn\\
	& \qquad \qquad \left. \rule{0pt}{13pt} \gamma_1 \left( \max\left\{ \beta_2(|x_2(0)|), \gamma_2 \left(\twonorm{w_2}\right)\right\} \right) \right\} \rule{1.0in}{0pt} \nn \\
	 & =  \max\left\{ \beta_1(|x_1(0)|), \rule{0pt}{13pt} \right. \nn\\
	 & \qquad \qquad \left. \rule{0pt}{13pt} \gamma_1 \circ \beta_2(|x_2(0)|), \gamma_1 \circ \gamma_2\left(\twonorm{w_2}\right) \right\}.
	 \label{eq:phi1_L2bnd}
\end{align}
For all $s \in \R_{\geq 0}$, define $\beta \in \Kinf$ by
 $\beta(s) \doteq \frac{1}{2} \max \left\{ \beta_1(s), \gamma_1 \circ \beta_2(s), \beta_2(s) \right\}$  
and $\gamma \in \Kinf$ by
$\gamma(s) \doteq \frac{1}{2} \max \left\{ \gamma_1 \circ \gamma_2(s), \gamma_2(s) \right\}$.
Combining \eqref{eq:z2_L2bnd} and \eqref{eq:phi1_L2bnd} yields
\begin{align}
	\twonorm{x}  & =  \twonorm{x_1} + \twonorm{x_2} \nn \\
	& \leq  \max \left\{ \beta(|x(0)|), \gamma\left(\twonorm{w_2}\right)\right\} \nn
\end{align}
and therefore the cascade connection of \eqref{eq:sys1}-\eqref{eq:sys2} with interconnection
$w_1 = x_2$ has the nonlinear $\L_2$-gain property.
$\halmos$

\subsection{Feedback Interconnection}
\label{sec:L2_feedback}
\begin{figure}[h]
\centering
\includegraphics[width=5.5cm]{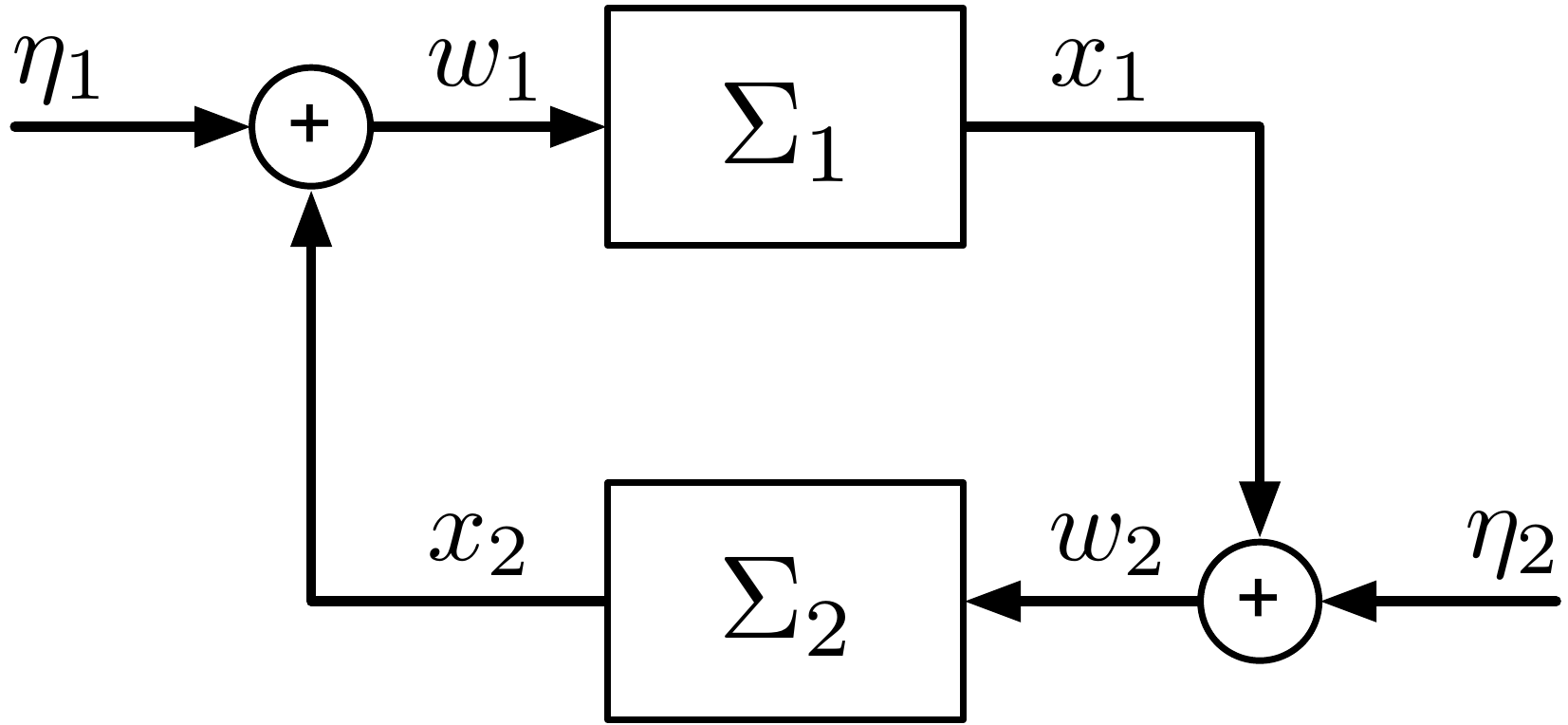}
\caption{Feedback Interconnection}
\label{fig:FbkDiag}
\end{figure}

We consider two feedback interconnections; one without external inputs 
(Figure~\ref{fig:FbkDiag} with $\eta_1 \equiv \eta_2 \equiv 0$) and one with external
inputs (Figure~\ref{fig:FbkDiag}, as shown).  We include the former as it has a much simpler
small-gain condition than the latter and gives rise to an interesting sufficient condition for the stability
of feedback interconnections of iISS systems, which we will discuss
 in Section~\ref{sec:iISS_connections} (see Theorem~\ref{thm:iISS_L2_noinput}).

In the diagram of Figure~\ref{fig:FbkDiag}, we obviously require that $\eta_i \in \W^{m_i}$, $i=1,2$, 
and, so that the input/output dimensions are consistent, we also require that $m_2 = n_1$ and $m_1 = n_2$.

\begin{thm}
	\label{thm:L2_fbk_noinputs}
	Suppose systems \eqref{eq:sys1} and \eqref{eq:sys2} satisfy the nonlinear \ltwo-gain
	bounds \eqref{eq:z1_L2bnd} and \eqref{eq:z2_L2bnd}, respectively, and the interconnection
	constraints $w_1 = x_2$ and $w_2 = x_1$.  
	 If the small-gain conditions
	\begin{equation}
		\label{eq:no_input_sgc}
		 {\rm Id} - \gamma_i \circ \gamma_j \in \Kinf
	\end{equation}
	are satisfied for $i,j = 1,2$, $i \neq j$, then the system is \ltwo-stable.
\end{thm}
{\it Proof:} For $i,j = 1,2$, $i\neq j$, using \eqref{eq:z1_L2bnd}, \eqref{eq:z2_L2bnd}, and the interconnection
constraints $w_i = x_j$,
we see that
\begin{align*}
	\lefteqn{\twonorm{x_i} \leq \max \left\{ \k_i(|x_i(0)|), \g_i\left(\twonorm{x_j}\right) \right\} } & & \nn \\
	& \leq  \max \left\{ \rule{0pt}{13pt} \k_i(|x_i(0)|), \right. \\
	& \qquad \qquad \left. \g_i \left( \max \left\{ \k_j(|x_j(0)|),  \g_j\left(
		\twonorm{x_i}\right) \right\} \right) \rule{0pt}{13pt} \right\}\nn \\
	& \leq  \max \left\{ \k_i(|x_i(0)|), \g_i  \circ \k_j(|x_j(0)|) \right\} \\
	& \qquad + \g_i \circ \g_j \left(
		\twonorm{x_i}\right).
\end{align*}
Therefore, if ${\rm Id} - \gamma_i  \circ \gamma_j \in \Kinf$, we can derive an
upper bound on $\twonorm{x} \doteq \twonorm{[x_1^T \, x_2^T]^T}$ 
that depends only on the initial condition
$x(0) \doteq [x_1^T(0) \, x_2^T(0)]^T$.  Specifically, let $\k_i \in \Kinf$ for $i,j=1,2$, $i \neq j$, be given by
\begin{multline*}
 	\bar\k_i(s) \doteq \max\left\{ \rule{0pt}{12pt}({\rm Id} - \g_i \circ \g_j ) \circ \k_i(s), \right. \\
	 \left. ({\rm Id} - \g_i \circ \g_j) \circ \g_i \circ \k_j(s) \rule{0pt}{12pt} \right\} , \quad \forall s \in \R_{\geq 0}. 
\end{multline*} 
  Furthermore, define $\k \in \Kinf$ by $\k \doteq \bar\k_1 + \bar\k_2$.  Then,
\begin{align}
	\twonorm{x}  & =  \twonorm{x_1} + \twonorm{x_2} \nn \\
	& \leq \bar\k_1(|x(0)|) + \bar\k_2(|x(0)|) 
	 \leq  \k(|x(0)|) \nn
\end{align}
demonstrating that the interconnected system is $\L_2$-stable. $\halmos$

\begin{thm}
\label{thm:small-gain-max}
	Suppose systems \eqref{eq:sys1} and \eqref{eq:sys2} satisfy the nonlinear \ltwo-gain
	bounds \eqref{eq:z1_L2bnd} and \eqref{eq:z2_L2bnd}, respectively, and the interconnection
	constraints $w_1 = x_2 + \eta_1$ and $w_2 = x_1 + \eta_2$.  Fix $\varepsilon \in \R_{>0}$.
	Let $\rho \in \Kinf$ be such that $(\rho - {\rm Id}) \in \Kinf$ and define
	\begin{equation}
	 \hat{\gamma}_i(s)  \doteq  \gamma_i \circ \rho((1+\varepsilon^2)s)
	 	\label{eq:hat_gamma}
	 \end{equation}
	for all $s \in \R_{\geq 0}$, $i=1,2$.  If the small-gain conditions
	\begin{equation}
		{\rm Id} - \hat{\gamma}_i \circ \hat{\gamma}_j  \in  \Kinf 
			\label{eq:sgc}
	\end{equation}
	are satisfied for $i,j=1,2$, $i \neq j$, then the interconnected system satisfies the nonlinear \ltwo-gain property
	from input $\eta = [\eta_1^T \ \eta_2^T]^T$ to state $x = [x_1^T \ x_2^T ]^T$.
\end{thm}

{\it Proof:}  Let $\mu \in \Kinf$ be given by $\mu \doteq \rho \circ (\rho - {\rm Id})^{-1}$.
We derive an upper bound on $\twonorm{x_i}$ using the nonlinear $\L_2$-gain property \eqref{eq:z1_L2bnd}, 
the interconnection condition $w_i = x_j + \eta_i$, Lemma~\ref{clm:Young}, and
Lemma~\ref{lem:weak_tri_ineq} as
\begin{align}
	\lefteqn{\twonorm{x_i}  \leq  \max \left\{ \k_i(|x_i(0)|), \g_i \left(\twonorm{w_i}\right) \right\} } & \nn\\
	 & =  \max \left\{ \k_i(|x_i(0)|), \g_i \left( \twonorm{x_j + \eta_i}\right) \right\} \nn \\
	& \leq  \max \left\{ \rule{0pt}{13pt} \k_i(|x_i(0)|), \right. \nn\\
	& \qquad \qquad \left. \g_i \left( (1+\varepsilon^2)\twonorm{x_j} + (1 + \tfrac{1}{\varepsilon^2})\twonorm{\eta_i} \right) \rule{0pt}{13pt} \right\} \nn \\
	& \leq  \max\left\{ \rule{0pt}{13pt} \k_i(|x_i(0)|), \hat{\g}_i \left(\twonorm{x_j}\right),  \right. \nn\\
	& \qquad \qquad \left. \g_i \circ \mu \left( (1+\tfrac{1}{\varepsilon^2})\twonorm{\eta_i}\right)
			\rule{0pt}{13pt} \right\} . \label{eq:int_1}
\end{align}
Repeating the same arguments, we can derive an upper bound on $\twonorm{x_j}$ that we then substitute into
\eqref{eq:int_1} to obtain
\begin{align*}
	\twonorm{x_i} & \leq  \max\left\{ \rule{0pt}{13pt} \k_i(|x_i(0)|), \hat{\g}_i \circ \k_j(|x_j(0)|), \right. \\
	& \qquad \qquad \left.  \g_i \circ \mu \left((1 + \tfrac{1}{\varepsilon^2})\twonorm{\eta_i} \right), \right. \nn\\
	&  \qquad \qquad \left.	\hat{\g}_i \circ \g_j \circ \mu\left((1+ \tfrac{1}{\varepsilon^2})\twonorm{\eta_j}\right) 
		\rule{0pt}{13pt} \right\} \\
	& \qquad +  \hat{\g}_i \circ \hat{\g}_j  \left( \twonorm{x_i}\right) .
\end{align*}
With the small gain condition \eqref{eq:sgc}, we see that we can upper bound $\twonorm{x_i}$ by terms
depending solely on initial conditions $x_1(0) \in \R^{n_1}$, $x_2(0) \in \R^{n_2}$ and inputs $\eta_1 \in \W^{m_1}$
and $\eta_2 \in \W^{m_2}$.  With the derived bounds on $\twonorm{x_1}$ and $\twonorm{x_2}$ we
may bound $\twonorm{x}$ as in the conclusion of the proofs of Proposition~\ref{prop:L2Cascade} and Theorem~\ref{thm:L2_fbk_noinputs} 
and we omit the details. $\halmos$

\begin{remark}
	We note that by choosing $\rho(s) = (1 + \varepsilon^2) s$ with $\varepsilon^2 <\!< 1$, 
	the small gain condition 
	\eqref{eq:sgc} approaches
	${\rm Id} - \gamma_1 \circ \gamma_2 ,  {\rm Id} - \gamma_2 \circ \gamma_1 \ \in \Kinf$.
	This is \eqref{eq:no_input_sgc} and the obvious analogue of the classical linear small-gain condition given by 
	$\gamma_1 \gamma_2 < 1$, with $\gamma_1, \gamma_2 \in \R_{>0}$.
	A further consequence of choosing $\varepsilon^2 <\! < 1$  is that the related functions or constants in Lemma~\ref{lem:weak_tri_ineq} 
	and Lemma~\ref{clm:Young} become large; i.e.,
	\[ \mu(s) \doteq \rho \circ (\rho - {\rm Id})^{-1}(s) = \frac{1+\varepsilon^2}{\varepsilon^2} s  \]
	and $1+\frac{1}{\varepsilon^2}$, respectively.
	As can be seen in the proof above, the function $\mu$ and
	the constant $1+\frac{1}{\varepsilon^2}$
	being large correspond to large
	bounds on external inputs.
	
	We observe that by using \eqref{eq:simple_tri_ineq} (i.e., $\rho(s) = 2s$) and $\varepsilon^2 = 1$
	the function in the small-gain condition reduces to
	$s - \hat\gamma_i  \circ \hat\gamma_j (s) = 
	s - \gamma_i \circ 4\gamma_j (4s) \in \Kinf$.
\end{remark}

Though we generally adhere to the maximum formulation of gain properties, rather than using the qualitatively
equivalent summation formulation, we state here a small-gain theorem for the $\L_2$-gain property given
by
\begin{equation}
	\label{eq:nonlinear_L2_plus}
	\twonorm{x_i} \leq \k_i(|x_i(0)|) + \g_i \left( \twonorm{w_i}\right).
\end{equation}
Despite the qualitative equivalence between \eqref{eq:nonlinear_L2_plus} and \eqref{eq:z1_L2bnd},
if one is interested in {\em quantitative} results, for example in looking to compute tight gain
bounds (e.g., \cite{ZhDo12-IJC}, \cite{ZhDo13-SCL}), the form of the small-gain condition below is useful.

\begin{thm}
\label{thm:small-gain-add}
	Suppose systems \eqref{eq:sys1} and \eqref{eq:sys2} satisfy the nonlinear \ltwo-gain
	bounds \eqref{eq:nonlinear_L2_plus} for $i=1,2$, respectively, and the interconnection
	constraints $w_1 = x_2 + \eta_1$ and $w_2 = x_1 + \eta_2$.  Fix $\varepsilon \in \R_{>0}$.
	Let $\rho \in \Kinf$ be such that $(\rho - {\rm Id}) \in \Kinf$ and define
	\begin{eqnarray}
	 \tilde{\gamma}_i (s) & \doteq & \gamma_i \circ \rho \circ \rho((1+\varepsilon^2)s) 
	 	\label{eq:tilde_gamma_add} \\
	 \hat{\gamma}_i(s) & \doteq & \gamma_i \circ \rho((1+\varepsilon^2)s)
	 	\label{eq:hat_gamma_add}
	 \end{eqnarray}
	for all $s \in \R_{\geq 0}$, $i=1,2$.  If the small-gain conditions
	\begin{equation}
		{\rm Id} - \tilde{\gamma}_i \circ \rho \circ \hat{\gamma}_j  \in  \Kinf 
	\end{equation}
	are satisfied for $i,j=1,2$, $i \neq j$, then the interconnected system satisfies the nonlinear \ltwo-gain property
	from input $\eta = [\eta_1^T \ \eta_2^T]^T$ to state $x = [x_1^T \ x_2^T ]^T$.
\end{thm}

\begin{remark}
	In order to avoid a proliferation of unnecessary notation in 
	Theorems~\ref{thm:small-gain-max} and \ref{thm:small-gain-add}, 
	we fixed a single constant
	$\varepsilon \in \R_{>0}$ and used a single function $\rho \in \Kinf$ such that
	$(\rho - {\rm Id}) \in \Kinf$.  In fact, each of the instances of these elements
	in \eqref{eq:hat_gamma}, \eqref{eq:tilde_gamma_add}, and \eqref{eq:hat_gamma_add}
	 may be chosen independently.
	For example, rather than a single $\varepsilon \in \R_{>0}$ in \eqref{eq:hat_gamma}, it is possible to
	choose two constants, say $\varepsilon_1, \varepsilon_2 \in \R_{>0}$.
	This follows from the fact that we could, in principle, fix a different $\varepsilon \in \R_{>0}$
	each time we apply Lemma~\ref{clm:Young} in the proof of Theorem~\ref{thm:small-gain-max}.
	A similar remark holds with regard to the two appearances of $\rho \in \Kinf$ corresponding
	to the two applications of Lemma~\ref{lem:weak_tri_ineq} in Theorem~\ref{thm:small-gain-max}, as well as for
	the constants $\varepsilon \in \R_{>0}$ and functions $\rho \in \Kinf$ in Theorem~\ref{thm:small-gain-add}.
\end{remark}



\section{Interconnections of (i)ISS Systems}
\label{sec:iISS_connections}
Cascade and feedback interconnections of ISS and iISS systems have been
extensively studied in the literature (see \cite{ChAn08-SCL, Ito10-TAC, ItJi09-TAC, RKW10-Aut}
and references therein).  In the case of iISS systems, it is known that iISS of the
individual subsystems alone is insufficient to guarantee desired properties such as
zero-input global asymptotic stability (0-GAS) or iISS of interconnected systems.
The results on interconnections of systems with nonlinear \ltwo-gain
in Section~\ref{sec:NonlinearL2gain} do not appear to require additional conditions
and, in light of the relationship between nonlinear \ltwo-gain and iISS
described in Theorem \ref{thm:iISSandNonlinear}, we now turn to the relationship
between interconnections of systems with nonlinear \ltwo-gain and known
results for the interconnection of ISS and iISS systems.
We first examine the feedback interconnection of ISS systems.
\begin{thm}
	\label{thm:ISS_L2_fbk}
	Suppose systems \eqref{eq:sys1} and \eqref{eq:sys2} are ISS with functions
	$\alpha_i,\beta_i,\sigma_i \in \Kinf$, $i=1,2$, as in \eqref{eq:ISS} and that
	the systems are connected in feedback with $w_1 = x_2$ and $w_2 = x_1$.
	Let $\bar{\gamma} = 1$ and let
	$T_i : \R^{n_i} \rightarrow \R^{n_i}$ and $S_i : \R^{m_i} \rightarrow \R^{m_i}$ be the 
	changes of coordinates from Theorem~\ref{thm:ISSandLinear} that yield new coordinates in which
	$\Sigma_1$ and $\Sigma_2$ satisfy linear \ltwo-gain bounds \eqref{eq:L2gain}
	with $\hat{\beta}_i \in \Kinf$ and $\bar\gamma = 1$.  If there exist $c_1, c_2 \in \R_{>0}$ such that,
	for $i,j = 1,2$, $i \neq j$,
	\begin{equation}
		|S_j(\zeta)|  \leq  \sqrt{c_i} |T_i(\zeta)|, \quad \forall \zeta \in \R^{n_i}, \label{eq:S2_bnd}
	\end{equation}
	and if $c_1 c_2 < 1$, then the feedback interconnection is $\alpha$-integrable.
\end{thm}

\noindent {\it Proof:} With the changes of coordinates for $\Sigma_1$ and $\Sigma_2$ that
yield  linear \ltwo-gain with $\bar\gamma=1$, 
\begin{eqnarray}
	\twonorm{T_i(x_i)} & \leq & \max\left\{ \hat\k_i (|T_i(x_i(0))|), \twonorm{S_i(w_i)} \right\} \nn 
\end{eqnarray}
for $i=1,2$.  Let $\psi_i \doteq T_i(x_i)$ and $\xi_i \doteq \psi_i(0) = T_i(x_i(0))$.  
Using the bounds \eqref{eq:S2_bnd},
 and the interconnection conditions, we obtain
\begin{align}
	\lefteqn{\twonorm{\psi_1}  \leq  \max \left\{ \hat\k_1(|T_1(x_1(0))|), \twonorm{S_1(x_2)} \right\} }& \nn\\
	& \leq  \max \left\{ \hat\k_1(|\xi_1|), c_2\twonorm{T_2(x_2)} \right\}   \nn \\
	& \leq  
	   \max \left\{ \hat\k_1(|\xi_1|), c_2 \hat\k_2(|\xi_2|), c_2 \twonorm{S_2(x_1)} \right\} \nn\\
	& \leq  \max\left\{ \hat\k_1(|\xi_1|), c_2 \hat\k_2(|\xi_2|) \right\} + c_2 c_1 \twonorm{\psi_1}.  \nn
\end{align}
Let $\xi \doteq [\xi_1 \, \xi_2]^T \in \R^{n_1+n_2}$ and $\k_i \in \Kinf$ for $i,j=1,2$, $i \neq j$, be defined by
\[ \k_i(s) \doteq \max\left\{ \hat\k_i(s), c_j \hat\k_j(s) \right\}, \quad \forall s \in \R_{\geq 0} . \]
Then
$\twonorm{\psi_1} \leq \frac{1}{1-c_1 c_2} \k_1(|\xi|)$
and a similar argument yields
$\twonorm{\psi_2} \leq \frac{1}{1-c_1 c_2} \k_2(|\xi|)$. 

Let $\alpha_\ell^i, \alpha_u^i \in \Kinf$ come from Lemma~\ref{lem:BoundChanges} applied
to the change of coordinates $T_i(\cdot)$, $i=1,2$, and let $\alpha \in \Kinf$ come from
Lemma~\ref{clm:LowerBoundKSum} applied to $\alpha_\ell^i \in \Kinf$.  Let $\rho \in \Kinf$ be such that
$\rho - {\rm Id} \in \Kinf$ and define $\mu \doteq \rho \circ (\rho - {\rm Id})^{-1} \in \Kinf$.
Define $\tilde\k, \k \in \Kinf$ by
\begin{eqnarray}
	\tilde\k(s) & \doteq & \frac{1}{1-c_1 c_2} \left( \k_1(s) + \k_2(s) \right) \nn \\
	\k(s) & \doteq & \max \left\{ \tilde\k \circ \rho \circ \alpha_u^1(s), \tilde\k \circ \mu \circ \alpha_u^2(s) \right\} \nn
\end{eqnarray}
for all $s \in \R_{\geq 0}$.
Then
\begin{align}
	\lefteqn{ \int_0^t \alpha(|x(\tau)|)d\tau  \leq  \int_0^t \alpha(|x_1(\tau)| + |x_2(\tau)|) d\tau } & \nn \\
	& \leq
	\int_0^t \alpha_\ell^1(|x_1(\tau)|)d\tau + \int_0^t \alpha_\ell^2(|x_2(\tau)|)d\tau \nn \\
	 & \leq    \twonorm{\psi_1} + \twonorm{\psi_2} \nn\\
	& \leq  \frac{1}{1-c_1 c_2} \left( \k_1(|\xi|) + \k_2(|\xi|) \right) 
	 =  \tilde\k(|\xi|) \nn \\
	 & \leq  \max \left\{ \tilde\k \circ \rho(|\xi_1|), \tilde\k \circ \mu(|\xi_2|) \right\} \nn \\
	& \leq  \max \left\{ \tilde\k \circ \rho \circ \alpha_u^1(|x_1(0)|), \tilde\k \circ \mu \circ \alpha_u^2(|x_2(0)|) \right\} \nn\\ 
	& =  \k(|x(0)|). \label{eq:alpha_int_bnd}
\end{align}
Therefore, the feedback interconnection is $\alpha$-integrable.
$\halmos$

The condition $c_1 c_2 < 1$ is analogous to the classical \ltwo \ small-gain
theorem \cite{vdSc00}.
While at first glance the above theorem may appear to provide a much simpler condition
to check than, for example, \cite[Theorem 2.1]{JTP94-MCSS}, finding the appropriate changes of coordinates
so that an arbitrary ISS system exhibits (linear) \ltwo-gain appears to be a difficult task.

When attempting to prove results on interconnected ISS systems
directly using the ISS estimates of Definition~\ref{def:ISS} it is sometimes necessary
to impose additional assumptions beyond those already known in the literature.
For example, the cascade interconnection of ISS systems is always ISS 
(e.g., \cite[Proposition 7.2]{Sont89-TAC}).  Attempting to prove this directly using the
ISS estimates of Definition~\ref{def:ISS} requires being able to
compare the state scaling, $\alpha_2 \in \Kinf$, of the
driving system with the input scaling, $\sigma_1 \in \Kinf$, of the driven system.
A similar assumption is required to prove a small-gain theorem using the
ISS estimates of Definition~\ref{def:ISS}.  Since such results are less general
than those available in the literature we do not present them here.


\subsection{Cascade Interconnections of iISS Systems}
\label{sec:Cascade_iISS}
In \cite{AAS02-SICON}, \cite{ChAn08-SCL}, and \cite{Ito10-TAC}, sufficient conditions are given guaranteeing
iISS of a cascade connection of iISS systems.  Generally, these conditions involve a relationship
between the decay rate of the driving system and the gain of the driven system
($\Sigma_2$ and $\Sigma_1$, respectively, of Figure~\ref{fig:CascadeDiag}).

The following sufficient condition for iISS of a cascade interconnection of iISS systems
makes use of the qualitative equivalence between iISS systems and those with nonlinear
\ltwo-gain as described in Theorem~\ref{thm:iISSandNonlinear}.
\begin{thm}
	\label{thm:cascade_iISS_L2}
	Suppose systems \eqref{eq:sys1}-\eqref{eq:sys2} are iISS with functions
	$\alpha_i,\beta_i,\gamma_i,\sigma_i \in \Kinf$, $i=1,2$, as in \eqref{eq:iISS} and that
	the systems are connected in cascade with $w_1 = x_2$.  Let 
	$T_i : \R^{n_i} \rightarrow \R^{n_i}$ and $S_i : \R^{m_i} \rightarrow \R^{m_i}$ be the 
	changes of coordinates from Theorem~\ref{thm:iISSandNonlinear} that yield new coordinates in which
	$\Sigma_1$ and $\Sigma_2$ satisfy nonlinear \ltwo-gain bounds \eqref{eq:nonlinearL2gain}
	with $\hat{\beta}_i, \hat{\gamma}_i \in \Kinf$.  If there exists $c \in \R_{>0}$ such that
	\begin{equation}
		\label{eq:casc_S1bnd}
		|S_1(\zeta)| \leq \sqrt{c} |T_2(\zeta)|, \quad \forall \zeta \in \R^{n_2}
	\end{equation}
	then the cascade interconnection is iISS.
\end{thm}

{\it Proof:}
Let $\psi_i \doteq T_i(x_i)$, $\xi_i \doteq \psi_i(0) = T_i(x_i(0))$, and $v_i \doteq S_i(w_i)$.
Using the bounds \eqref{eq:z1_L2bnd} and \eqref{eq:z2_L2bnd} with $\hat\k_i, \hat\gamma_i \in \Kinf$,
the interconnection condition $w_1 = x_2$, the bound \eqref{eq:casc_S1bnd}, and 
Lemma~\ref{lem:weak_tri_ineq},
we obtain
\begin{align}
	\lefteqn{\twonorm{\psi_1}  
	 \leq  \max \left\{ \hat\k_1(|\xi_1|), \hat\g_1\left(\twonorm{S_1(x_2)}\right) \right\} } & \nn\\
	 & \leq  \max \left\{ \hat\k_1(|\xi_1|), \hat\g_1\left(c\twonorm{T_2(x_2)}\right) \right\}   \nn\\
	& \leq   \max \left\{ \hat\k_1(|\xi_1|), \hat\g_1\left(c\hat\k_2(|\xi_2|)\right),  \hat\g_1 \left(c \hat\g_2\left(\twonorm{v_2}\right)\right) \right\}. \rule{2.0in}{0pt} \nn 
\end{align}
For $i=1,2$, let $\alpha_\ell^i \in \Kinf$ come from Lemma~\ref{lem:BoundChanges} 
applied to $T_i(\cdot)$ and let $\tilde\alpha \in \Kinf$ come
from Lemma~\ref{clm:LowerBoundKSum} applied to $\alpha_\ell^i \in \Kinf$.  This yields the following
bound:
\begin{align}
	\lefteqn{\int_0^t \tilde\alpha (|x(\tau)|) d\tau 
	 \leq  \int_0^t \left(\alpha_\ell^1(|x_1(\tau)|)
		+ \alpha_\ell^2(|x_2(\tau)|)\right) d\tau } & \nn\\
	& \leq  \twonorm{\psi_1} + \twonorm{\psi_2}  \nn \\
	& \leq  2 \max \left\{ \rule{0pt}{14pt} \hat\k_1(|\xi_1|), \hat\g_1\left(c\hat\k_2(|\xi_2|)\right), \right. \nn\\
	& \qquad \qquad \quad \left. \hat\g_1 
		\left(c \hat\g_2\left(\twonorm{v_2}\right)\right), 
	 \hat\k_2(|\xi_2|), \right. \nn\\
	 & \qquad \qquad \quad \left.  \hat\g_2\left(\twonorm{v_2}\right) \right\} .
	\label{eq:casc_bnd}
%
\end{align}
For all $s \in \R_{\geq 0}$, define $\tilde\k \in \Kinf$ by
	$\tilde\k(s) \doteq 2 \max \left\{ \hat\k_1 \circ \alpha_u^1(s),
		\hat\g_1(c \hat\k_2 \circ \alpha_u^2(s)), \hat\k_2 \circ \alpha_u^2(s) \right\}$
and $\tilde\g \in \Kinf$ as
	$\tilde\g(s) \doteq 2 \max \left\{ \hat\g_1(c\hat\g_2(s)), \hat\g_2(s) \right\}$.
Finally, Lemma~\ref{lem:BoundChanges} allows us to upper bound $|S_2(w_2)|$
by $\tilde\sigma^{1/2} \in \Kinf$.  With the above definitions, \eqref{eq:casc_bnd} becomes
\[ \int_0^t \tilde\alpha (|x(\tau)|) d\tau \leq \max \left\{ \tilde\k(|x(0)|), \tilde\g \left(\int_0^t 
	\tilde\sigma(|w_2(\tau)|) d\tau \right) \right\} \]
which demonstrates that the cascade is iISS from input $w_2$ to state $x$. $\halmos$


From Theorem~\ref{thm:cascade_iISS_L2}
we note that while iISS and nonlinear \ltwo-gain are qualitatively equivalent 
properties (Theorem~\ref{thm:iISSandNonlinear}), there is no contradiction between the extra conditions required to guarantee
iISS of cascaded systems in \cite{AAS02-SICON, ChAn08-SCL} and the lack of such extra conditions in Proposition~\ref{prop:L2Cascade}
for the cascade of systems with nonlinear $\L_2$-gain.  
In fact, the result of Theorem~\ref{thm:cascade_iISS_L2} is
similar to the results of \cite{AAS02-SICON, ChAn08-SCL} in that the sufficient condition 
requires a relationship between
the state change of coordinates of the driving system and the input change of coordinates
of the driven system.

The qualitatively equivalent definition of iISS in Definition~\ref{def:iISS} (\cite{ASW00-DC})
gives rise to the following similar sufficient condition for iISS of cascaded iISS systems.
\begin{thm}
	\label{thm:eq_cascade_iISS}
	Suppose systems \eqref{eq:sys1}-\eqref{eq:sys2} are iISS with functions
	$\alpha_i,\beta_i,\gamma_i,\sigma_i \in \Kinf$, $i=1,2$, as in \eqref{eq:iISS} and that
	the systems are connected in cascade with $w_1 = x_2$.  If there exists a $c \in \R_{>0}$ so
	that
	\begin{equation}
		\label{eq:cascade_iISS_cdn}
		\sigma_1(s) \leq c \alpha_2(s), \quad \forall s \in \R_{\geq 0},
	\end{equation}
	then the cascade interconnection is iISS.
\end{thm}

We note that the simplicity of the condition \eqref{eq:cascade_iISS_cdn} and the following
proof, as compared with the results of \cite{AAS02-SICON} or \cite{ChAn08-SCL},
stems from the fact that the iISS property defined by \eqref{eq:iISS} treats the input
and the state in the same manner; i.e., \eqref{eq:iISS} is an integral-to-integral
estimate.  By contrast, \eqref{eq:iISS_original} does not treat the input and state in the
same manner and consequently relating the input of the driven system to the state of the
driving system requires more involved arguments.

{\it Proof:}
The iISS estimate for system \eqref{eq:sys2} is given by
\begin{multline}
	\int_0^t \alpha_2(|x_2(\tau)|) \\
	\leq \max \left\{ \beta_2(|x_2(0)|), 
	 \gamma_2 \left(\int_0^t \sigma_2(|w_2(\tau)|)d\tau\right) \right\} . \label{eq:Sigma2_bnd}
\end{multline}
With the iISS estimate \eqref{eq:iISS} for systems \eqref{eq:sys1}-\eqref{eq:sys2},
the interconnection condition $w_1 = x_2$, and the condition \eqref{eq:cascade_iISS_cdn}, 
 the following
calculation is straightforward:
\begin{multline}
	\int_0^t \alpha_1(|x_1(\tau)|)d\tau \leq \\
	   \max \left\{ \rule{0pt}{16pt} \beta_1(|x_1(0)|), \gamma_1 \left(c\beta_2(|x_2(0)|)\right), \right. \\
		\left. \rule{0pt}{16pt} \gamma_1 \left( c \gamma_2 \left(
		\int_0^t \sigma_2(|w_2(\tau)|)d\tau\right)\right) \right\}. \label{eq:Sigma1_bnd}
\end{multline}

Let $\alpha \in \Kinf$ come from the application of Lemma~\ref{clm:LowerBoundKSum}
to $\alpha_1,\alpha_2 \in \Kinf$ so that
\eqref{eq:sum_bound} together with the 
bounds \eqref{eq:Sigma2_bnd} and \eqref{eq:Sigma1_bnd} implies
\begin{equation}
	\int_0^t \alpha(|x(t)|)d\tau \leq \max \left\{ \beta(|x(0)|), 
	 \gamma\left(\int_0^t \sigma_2(|w_2(\tau)|)d\tau\right) \right\} \nn
\end{equation}
where, for all $s \in \R_{\geq 0}$, $\k,\g \in \Kinf$ are given by
$\k(s)  \doteq  2\max\left\{ \k_1(s), \g_1(c\k_2(s)),\k_2(s) \right\}$ and
$\g(s)  \doteq  2\max \left\{ \g_1(c\g_2(s)), \g_2(s) \right\}$.  Therefore,
 the cascade system is iISS from input $w_2$ to state
$x = [x_1 \, x_2]^T$. $\halmos$


\subsection{Feedback Interconnections of iISS Systems}
Similar to the extra conditions required for cascades of iISS systems to be iISS,
sufficient conditions for iISS of feedback interconnections have been shown
to require more than iISS of the subsystems and a small-gain condition
(see \cite{Ito06-TAC}, \cite{ItJi09-TAC}).  We now turn to the relationship between
iISS systems, coordinates in which these systems satisfy the nonlinear \ltwo-gain property, and the
small-gain result of Theorem~\ref{thm:small-gain-add}.

\begin{thm}
	\label{thm:iISS_L2_noinput}
	Suppose systems \eqref{eq:sys1}-\eqref{eq:sys2} are iISS with functions
	$\alpha_i,\beta_i,\gamma_i,\sigma_i \in \Kinf$, $i=1,2$, as in \eqref{eq:iISS} and that
	the systems are connected in feedback with $w_1 = x_2$ and $w_2 = x_1$.  Let 
	$T_i : \R^{n_i} \rightarrow \R^{n_i}$ and $S_i : \R^{m_i} \rightarrow \R^{m_i}$ be the 
	changes of coordinates from Theorem~\ref{thm:iISSandNonlinear} that yield new coordinates in which
	$\Sigma_1$ and $\Sigma_2$ satisfy nonlinear \ltwo-gain bounds \eqref{eq:nonlinearL2gain}
	with $\hat{\beta}_i, \hat{\gamma}_i \in \Kinf$.  If there exist $c_1,c_2 \in \R_{>0}$ such that,
	for $i,j=1,2$, $i \neq j$,
	\begin{equation}
		|S_j(\zeta)|  \leq  \sqrt{c_i} |T_i(\zeta)|, \quad \forall \zeta \in \R^{n_i}, \label{eq:s2_bnd}
	\end{equation}
	and if the small-gain conditions
	\begin{equation}
		\label{eq:iISS_L2_sgc}
		{\rm Id} - \hat{\gamma}_i \left( c_j \hat{\gamma}_j(c_i \cdot) \right)  \in  \Kinf 
	\end{equation}
	hold, then the feedback interconnection is $\alpha$-integrable.
\end{thm}

\begin{remark}
We observe that, with the exception of the constants $c_1$ and $c_2$, the small-gain condition in
Theorem~\ref{thm:iISS_L2_noinput} is the same as that of Theorem~\ref{thm:L2_fbk_noinputs}.
However, to obtain a stability result for general iISS systems requires the additional 
conditions on state and input changes of coordinates given in 
\eqref{eq:s2_bnd}.
\end{remark}

 {\it Proof:} With the changes of coordinates for $\Sigma_1$ and $\Sigma_2$ that
yield  nonlinear \ltwo-gain with $\hat\k_i,\hat\g_i \in \Kinf$, we have
\begin{multline}
	\twonorm{T_i(x_i)}  \\
	\leq  \max\left\{ \hat\k_i (|T_i(x_i(0))|), \hat\g_i\left(\twonorm{S_i(w_i)}\right) \right\} \nn 
\end{multline}
for $i=1,2$.  Let $\psi_i \doteq T_i(x_i)$ and $\xi_i \doteq \psi_i(0) = T_i(x_i(0))$.  Using the bounds \eqref{eq:s2_bnd},
 and the interconnection conditions we obtain
\begin{align}
	\twonorm{\psi_1}  & \leq  \max \left\{ \hat\k_1(|\xi_1)|), \hat\g_1\left(\twonorm{S_1(x_2)} \right)\right\} 
	  \nn\\
	& \leq  \max \left\{ \hat\k_1(|\xi_1|), \hat\g_1\left(c_2\twonorm{T_2(x_2)} \right) \right\}   \nn \\
	& \leq  \max\left\{ \hat\k_1(|\xi_1|), \hat\g_1\left(c_2 \hat\k_2(|\xi_2|)\right), \rule{0pt}{12pt} \right. \nn\\
	& \qquad	\qquad \left. \rule{0pt}{12pt} \hat\g_1\left(c_2 \hat\g_2 \left( \twonorm{S_2(w_2)} \right) \right) \right\}  \nn \\
	& \leq  \max\left\{ \hat\k_1(|\xi_1|), \hat\g_1(c_2 \hat\k_2(|\xi_2|)) \right\} \nn\\
	& \qquad + \hat\g_1 \left(
			c_2 \hat\g_2 \left( c_1 \twonorm{\psi_1} \right) \right). \nn
\end{align}
Then the small-gain condition \eqref{eq:iISS_L2_sgc} yields an upper bound on 
$\twonorm{\psi_1}$ in terms of the initial condition $\xi \in \R^{n_1+n_2}$.  A similar argument
yields a similar bound for $\twonorm{\psi_2}$.  From here we follow the argument in 
\eqref{eq:alpha_int_bnd} to obtain that the feedback interconnection is $\alpha$-integrable. $\halmos$

When we additionally allow external inputs we obtain the following result.
\begin{thm}
	\label{thm:iISS_L2_fbk}
	Suppose systems \eqref{eq:sys1}-\eqref{eq:sys2} are iISS with functions
	$\alpha_i,\beta_i,\gamma_i,\sigma_i \in \Kinf$, $i=1,2$, as in \eqref{eq:iISS} and that
	the systems are connected in feedback with $w_1 = x_2 + \eta_1$ and $w_2 = x_1 + \eta_2$.  Let 
	$T_i : \R^{n_i} \rightarrow \R^{n_i}$ and $S_i : \R^{m_i} \rightarrow \R^{m_i}$, $i=1,2$, be the 
	changes of coordinates from Theorem~\ref{thm:iISSandNonlinear} that yield new coordinates in which
	$\Sigma_1$ and $\Sigma_2$ satisfy nonlinear \ltwo-gain bounds \eqref{eq:nonlinearL2gain}
	with $\hat{\beta}_i, \hat{\gamma}_i \in \Kinf$.  Fix $\varepsilon \in \R_{>0}$ and let $\rho \in \Kinf$ 
	be such that $\rho - {\rm Id} \in \Kinf$.  Suppose there exist $c_{S_i},c_{T_i} \in \R_{>0}$,
	$i=1,2$, such that
	\begin{eqnarray}
		|S_i(\zeta)|  \leq  \sqrt{c_{S_i}} |\zeta| & , & \quad \forall \zeta \in \R^{m_i}, \label{eq:S_sectors}\\
		|\zeta|  \leq  \sqrt{c_{T_i}} |T_i(\zeta)|& , & \quad \forall \zeta \in \R^{n_i}, \label{eq:T_sectors}
	\end{eqnarray}
	and define $\tilde\gamma_i \in \Kinf$, $i,j=1,2$, $i\neq j$, by
	\begin{equation}
		\tilde\gamma_i (s) \doteq \hat\gamma_i \circ \rho \left(c_{S_i} c_{T_j} (1+\varepsilon^2)s\right), \ \forall
		s \in \R_{\geq 0}.
	\end{equation}
	If the small-gain conditions
	\begin{equation}
		 {\rm Id} - \tilde\gamma_i \circ \tilde\gamma_j \in \Kinf 
	\end{equation}
	hold for $i=1,2$, $i\neq j$, then the feedback interconnection is iISS.
\end{thm}

The proof is similar to that of Theorem~\ref{thm:iISS_L2_noinput} with two essential differences.
The first is in the need to appeal to Lemma~\ref{lem:weak_tri_ineq} to overbound
$\Kinf$ functions of sums.  The second difference comes from the need to use
Lemma~\ref{clm:Young} and conditions \eqref{eq:S_sectors} and \eqref{eq:T_sectors}
to overbound \ltwo-norms of sums.  In particular, the sector bounds \eqref{eq:S_sectors}
and \eqref{eq:T_sectors} are used to derive bounds in the following way:
\begin{align}
	\lefteqn{ \twonorm{S_i(w_i)}  =  \twonorm{S_i(x_j + \eta_i)} } & \nn \\
	 & \leq  c_{S_i} \twonorm{x_j + \eta_i} \nn \\
	& \leq  c_{S_i} (1+ \varepsilon^2)\twonorm{x_j} + c_{S_i} (1 + \tfrac{1}{\varepsilon^2})
		\twonorm{\eta_i} \nn \\
	& \leq   c_{S_i} (1+ \varepsilon^2)c_{T_j}\twonorm{T_j(x_j)} + c_{S_i} (1 + \tfrac{1}{\varepsilon^2})
		\twonorm{\eta_i} \nn .
\end{align}
We then appeal to the fact that, in the coordinates defined by the change of coordinates
$T_j(\cdot)$, the system has the nonlinear \ltwo-gain property.  With this calculation, the proof
then closely follows that of Theorem~\ref{thm:iISS_L2_noinput} and we omit further details.

Similar to Theorem~\ref{thm:eq_cascade_iISS}, the qualitatively equivalent definition of iISS
in Definition~\ref{def:iISS} (\cite{ASW00-DC}) yields the following novel
 sufficient condition for iISS of the feedback interconnection of iISS systems.
\begin{thm}
	\label{thm:iISS_fbk}
	Suppose systems \eqref{eq:sys1}-\eqref{eq:sys2} are iISS with functions
	$\alpha_i,\beta_i,\gamma_i,\sigma_i \in \Kinf$, $i=1,2$, as in \eqref{eq:iISS} and that
	the systems are connected in feedback with $w_1 = x_2 + \eta_1$ and $w_2 = x_1 + \eta_2$.
	If there exist $c_i \in \R_{>0}$ and 
	$\rho_i \in \Kinf$ with $\rho_i - {\rm Id} \in \Kinf$ so that
	\begin{equation}
		\label{eq:iISS_fbk_constraint}
		\sigma_i \circ \rho_i (s) \leq c_j \alpha_j(s), \quad \forall s \in \R_{\geq 0}
	\end{equation}
	for $i,j=1,2$, $i\neq j$, and if there exists $\rho \in \Kinf$ with $\rho - {\rm Id} \in \Kinf$ so that
	\begin{equation}
		{\rm Id} - \gamma_i  \circ \rho \circ c_j \gamma_j  (c_i \cdot )\in \Kinf, 
	\end{equation}
	then the feedback interconnection is iISS.
\end{thm}

The proof of Theorem~\ref{thm:iISS_fbk} involves deriving several lengthy but entirely straightforward
upper bounds similar to the proofs of Theorems \ref{thm:L2_fbk_noinputs}, \ref{thm:small-gain-max}, and
\ref{thm:eq_cascade_iISS}. We omit the details.

\section{Conclusion}
\label{sec:Conc}
In this paper we have clarified the relationship between various ISS-type and $\L_2$-type stability properties.
In particular, we have demonstrated the qualitative equivalence between $\L_2$-stability and $\alpha$-integrability,
between linear $\L_2$-gain and ISS, and between nonlinear $\L_2$-gain and integral ISS.  Demonstrating these qualitative 
equivalences is done by considering nonlinear changes of coordinates.

We further presented several new sufficient conditions for stability
 of systems connected in cascade or feedback.
These conditions are derived using various qualitatively equivalent versions of the desired stability 
properties.  In particular, we have clarified the relationship between 
cascade and feedback stability results for systems with nonlinear \ltwo-gain and results
in the literature for cascade and feedback stability results for iISS systems, where the 
latter systems are known to require extra conditions to guarantee the desired stability results.

It appears unlikely that the results of Theorems \ref{thm:cascade_iISS_L2}, 
\ref{thm:iISS_L2_noinput}, and \ref{thm:iISS_L2_fbk} will be useful in
the sense of providing easily checkable conditions for stability of cascade or feedback
interconnected systems due to the fact that finding the appropriate changes of coordinates
would seem to be a challenging task.  However, these theorems serve to illustrate that there
is no contradiction between the stability results for systems already in coordinates such that
the nonlinear \ltwo-gain property holds, such as Proposition \ref{prop:L2Cascade} and 
Theorem \ref{thm:small-gain-add}, and
the results available in the literature on interconnections of iISS systems.

\section*{Acknowledgements}
The authors would like to thank Hiroshi Ito, Fabian Wirth, and Huan Zhang for helpful discussions
on this work.

\bibliographystyle{abbrvurl}
\bibliography{Control2011,ControlBooks}

\end{document}